\documentclass[a4paper,12pt]{article}

\usepackage[english]{babel}
\usepackage{amsmath}
\usepackage{amssymb}
\usepackage{epsfig}
\usepackage{theorem}
\usepackage{enumerate}
\usepackage{hyperref}
\usepackage{calc}
\usepackage{amsfonts}
\usepackage{ifthen}

\newcounter{count}[section]
\renewcommand{\thecount}{\thesection.\arabic{count}}

\newcommand{\equa}[1]{\addtocounter{count}{1}\begin{equation}#1\end{equation}}
\newcommand{\eqna}[1]{\addtocounter{count}{1}\begin{eqnarray}#1\end{eqnarray}}
\newcommand{\alig}[1]{\addtocounter{count}{1}\begin{align}#1\end{align}}

\newcounter{tictac1}
\newcounter{tictac2}
\newcounter{tictac3}
\newcounter{tictac4}
\newenvironment{fleuve1}{
   \begin{list}{\rm{$\textbf{\arabic{tictac1})}$} }{\usecounter{tictac1}
\leftmargin 1cm\labelwidth 2em}}{\end{list}}
\newenvironment{fleuvea}{
   \begin{list}{\rm{$\textbf{(\alph{tictac2})}$} }{\usecounter{tictac2}
\leftmargin 1cm\labelwidth 2em}}{\end{list}}
\newenvironment{fleuveI1}{
   \begin{list}{\rm{$\textbf{(I\arabic{tictac3})}$} }{\usecounter{tictac3}
\leftmargin 1cm\labelwidth 2em}}{\end{list}}
\newenvironment{fleuvei}{
   \begin{list}{\rm{$\textbf{(\roman{tictac4})}$} }{\usecounter{tictac4}
\leftmargin 1cm\labelwidth 2em}}{\end{list}}

\newlength{\Longueur}
\newenvironment{remas1}[1]{\smallskip\par\noindent\textbf{Remarks}\vspace{-.25cm}
\ifthenelse{#1=0}{\settowidth{\Longueur}{00.00}}{\addtocounter{count}{#1}
\settowidth{\Longueur}{\thecount}\addtocounter{count}{-#1}}
\begin{list}{\addtocounter{count}{1}\thecount}{\setlength{\labelsep}{.2cm}\setlength{\labelwidth}
{\the\Longueur}\setlength{\leftmargin}{\labelwidth+\labelsep}}}{\end{list}}

\theoremstyle{break}
\newtheorem{theo}[count]{Theorem}

\newtheorem{prop}[count]{Proposition}

\newtheorem{defi}[count]{Definition}
\newtheorem{lemm}[count]{Lemma}
\theorembodyfont{\rmfamily}
\newtheorem{rema}[count]{Remark}

\newcommand{\var}[1]{{\rm var}_{#1}}

\newcommand{\Spr}{{\rm Spr}}

\newcommand{\Id}{{\rm Id}}
\newcommand{\reff}[1]{(\ref{#1})}

\newcommand{\proof}{\noindent\textbf{Proof}\quad}

\def\1{\,\rlap{\mbox{\small\rm 1}}\kern.15em 1}
\def\ind#1{\1_{#1}}
\def\build#1_#2^#3{\mathrel{\mathop{\kern 0pt#1}\limits_{#2}^{#3}}}
\def\tend#1#2#3{\build\hbox to 12mm{\rightarrowfill}_{#1\rightarrow
#2}^{#3}}

{}
\def\converge#1#2#3{\build\hbox to
15mm{\rightarrowfill}_{\hbox{\scriptsize #3}}^{#1\rightarrow #2}}
\def\converg#1#2#3{\build\hbox to
15mm{\rightarrowfill}_{\hbox{\scriptsize #3}}^{#1\uparrow #2}}
\def\embf#1{\emph{\bf #1}}

\title{\bf{Chains with complete connections and one-dimensional Gibbs measures}}
\author{Roberto Fern\'andez\thanks{roberto.fernandez@univ-rouen.fr}\\
Gr\'egory Maillard\thanks{gregory.maillard@univ-rouen.fr}\\
\small{{\it Laboratoire de Math\'ematiques Rapha{\"e}l Salem}}\\ \small{\it{UMR 6085 CNRS-Universit\'e de Rouen}}\\
\small{\it{Site Colbert F-76821 Mont Saint Aignan, France}}}
\begin{document}
%\selectlanguage{american}
\maketitle
{\begin{abstract}
We discuss the relationship between discrete-time processes (chains) and one-dimensional Gibbs measures.  We consider
finite-alphabet (finite-spin) systems, possibly with a grammar (exclusion rule).  We establish conditions for a
stochastic process to define a Gibbs measure and vice versa.  Our conditions generalize well known equivalence results
between ergodic Markov chains and fields, as well as the known Gibbsian character of processes with exponential
continuity rate.  Our arguments are purely probabilistic; they are based on the study of regular systems of conditional
probabilities (specifications).  Furthermore, we discuss the equivalence of uniqueness criteria for chains and fields
and we establish bounds for the continuity rates of the respective systems of finite-volume conditional probabilities.
As an auxiliary result we prove a (re)construction theorem for specifications starting from single-site conditioning,
which applies in a more general setting (general spin space, specifications not necessarily Gibbsian).
\end{abstract}} \thispagestyle{empty}
%\tableofcontents \thispagestyle{empty} \newpage

\section{Introduction}
One dimensional systems are simultaneously the object of the theory of stochastic processes and the theory of Gibbs
measures. The complementarity of both approaches has yet to be fully exploited. Stochastic processes are defined on the
basis of transition probabilities.  A consistent chain is one for which these probabilities are a realization of the
single-site conditional probabilities given the past.  A Gibbs measure is defined in terms of specifications, which
determine its finite-volume conditional probabilities given the exterior of the volume.  In one dimension this implies
conditioning both the past and the future.  In this paper we study conditions under which a stochastic process defines,
in fact, a Gibbs measure and, in the opposite direction, when a Gibbs measures can be seen as a stochastic process.

This type of questions has been completely elucidated for Markov processes and fields.  See, for instance, Chapter 11
of the treatise by Georgii (1988)\nocite{geo88}. The equivalence, however, is obtained by eigenvalue-eigenvector
considerations which are not readily applicable to non-Markovian processes.  The Gibbsian character of processes with
exponentially decreasing continuity rate is also known. It follows from Bowen's characterization of Gibbs measures
(Theorem 5.2.4 in Keller, 1998\nocite{kel98}, for instance).  No result seems to be available on the opposite
direction, namely on the characterization of a one-dimensional Gibbs measure for an exponentially summable interaction
as a stochastic process.

In our paper we present both a generalization and an alternative to this previous work.  We directly establish
consistency-preserving maps between specifications and transition probabilities.  More precisely, these applications
are between specifications and their analogous for stochastic processes, which we call \emph{left-interval
specifications} (LIS).  The description in terms of LIS is equivalent to that in terms of transition probabilities, but
it offers a setting that mirrors the statistical mechanical setting of Gibbs measures.  In fact, the use of LIS allows
us to ``import'', in a painless manner, concepts and results from statistical mechanics into the theory of stochastic
processes.  This will be further exploited in a companion paper (Fern\'andez and Maillard, 2003).

We consider systems with a finite alphabet, possibly with a grammar, that is, with exclusion rules such that the
non-excluded configurations form a compact set.  We do not assume translation invariance either of the kernels or of
the consistent measures.  The main limitation of our results is that, in order to insure that the necessary limits are
uniquely defined, specifications and processes are required to satisfy a strong uniqueness condition called
\emph{hereditary uniqueness condition} (HUC).  A second property, called \emph{good future} (GF) is demanded for
stochastic processes to guarantee some control of the conditioning with respect to the future.  HUC is verified, for
instance, by specifications satisfying Dobrushin and boundary-uniformity criteria (reviewed below). Both GF and HUC are
satisfied by a large family of processes, for instance by the chains with summable variations studied by Harris
(1955)\nocite{har55}, Ledrappier (1974)\nocite{led74}, Walters (1975)\nocite{wal75},Lalley (1986)\nocite{lal86}, Berbee
(1987)\nocite{ber87}, Bressaud et al (1999)\nocite{bfg2}, \ldots.

Our results show that under these hypotheses there exist: (i) a map that to each LIS associates a specification such
that the process consistent with the former is a Gibbs measure consistent with the latter (Theorem \ref{th1}), and (ii)
a map that to each specification associates a LIS such that the Gibbs measure consistent with the former is a process
consistent with the latter (Theorem \ref{th2}). If domain and image match, these maps are inverses of each other. This
happens, in particular, in the case of exponentially decreasing continuity rates (Theorem \ref{th4}).  As part of the
proofs, we obtain estimates linking the continuity rates of LIS and specifications related by these maps (Theorem
\ref{th5}).  We also show that the validity of the Dobrushin and boundary-uniformity criteria for the specification
implies the validity of analogous criteria for the associated stochastic process (Theorem \ref{th3}). Finally, in
Appendix \ref{A2} we show that a system of single-site normalized kernels, satisfying order-consistency and boundedness
properties with respect to an a-priori measure, can be extended, in a unique manner, to a full specification.  This
generalizes the reconstruction Theorem 1.33 in Georgii (1988)\nocite{geo88}.  As this theorem may be of independent
interest, we have stated it in rather general terms, for arbitrary spin space and non-necessarily Gibbsian kernels
(Theorem \ref{thap}).

\section{Notation and preliminary definitions}\label{sect1}

We consider a finite alphabet $\mathcal{A}$ endowed with the discrete topology and $\sigma$-algebra, and $\Omega$ a
compact subset of $\mathcal{A}^{\mathbb{Z}}$. The space $\Omega$ is endowed with the projection $\mathcal{F}$ of the
product $\sigma$-algebra associated to $\mathcal{A}^{\mathbb{Z}}$. The space $\Omega$ represents admissible ``letter
configurations", where the admissibility is defined, for instance by some exclusion rule as in Ruelle
(1978)\nocite{rue78} or by a ``grammar" (subshift of finite type) as in Walters (1975)\nocite{wal75}.  For each
$\Lambda\subset\mathbb{Z}$, and each configuration $\sigma\in\mathcal{A}^{\mathbb{Z}}$ we denote $\sigma_\Lambda$ its
projection on $\Lambda$, namely the family $(\sigma_i)_{i\in\Lambda}\in \mathcal{A}^{\Lambda}$.  We denote
\equa{
  \label{eq:1} \Omega_\Lambda \;\triangleq\;
\Bigl\{\sigma_\Lambda\in\mathcal{A}^{\Lambda}: \exists\, \omega\in
\Omega \hbox{ with } \omega_\Lambda=\sigma_\Lambda \Bigr\} \;.
}
We denote $\mathcal{F}_{\Lambda}$ the corresponding
sub-$\sigma$-algebra of $\mathcal{F}$.
%, generated by the cylinders
%\begin{equation}
%C(\sigma_\Delta) \;\triangleq\; \Bigl\{ \omega\in\Omega :
%\omega_\Delta=\sigma_\Delta \Bigr\}
%\end{equation}
%for $\Delta\subset\Lambda$ and $\sigma_\Delta\in\Omega_\Delta$.
When $\Lambda$ is an interval, $\Lambda=[k,n]$ with $-\infty\le k\le
n\le +\infty$, we shall use the ``sequence'' notation:
$\omega_{k}^{n}\triangleq\omega_{[k,n]}=\omega_{k}, \ldots ,
\omega_{n}$, $\Omega_k^n\triangleq\Omega_{[k,n]}$, etc.  The notation
$\omega_\Lambda\, \sigma_\Delta$, where $\Lambda\cap\Delta=\emptyset$,
indicates the configuration on $\Lambda\cup\Delta$ coinciding with
$\omega_i$ for $i\in\Lambda$ and with $\sigma_i$ for $i\in\Delta$.  In
particular, $\omega_{k}^{n} \, \sigma_{n+1}^{m} = \omega_{k},
\ldots, \omega_{n}, \sigma_{n+1}, \ldots , \sigma_{m}$.
For $\omega,\sigma \in \mathcal{A}^{\mathbb{Z}}$, we note
\equa{ \label{ccc5} \sigma \stackrel{\neq j}{=} \omega
\quad \Longleftrightarrow \quad \sigma_{i}=\omega_{i}\;, \; \forall \,
i\neq j}
(``$\sigma$ equal to $\omega$ off $j$'').

We denote $\mathcal{S}$ the
set of finite subsets of $\mathbb{Z}$ and $\mathcal{S}_{b}$ the set of
finite intervals of $\mathbb{Z}$.
% and $\mathcal{S}_{b}^{\infty}$ the set of intervals of
% $\mathbb{Z}$.
For every $\Lambda \in \mathcal{S}_{b}$ we denote $l_{\Lambda}
\triangleq \min \Lambda$ and $m_{\Lambda} \triangleq \max \Lambda$,
$\Lambda_{-} = ]-\infty,l_\Lambda - 1]$, $\Lambda_{+} =[m_\Lambda +1,
+\infty[$ and $\Lambda_{+}^{(k)}=[m_\Lambda+1, m_\Lambda +k]$ for all
$k \in \mathbb{N}^{*}$.  The expression $\lim_{\Lambda \uparrow V}$
will be used in two senses. For kernels associated to a LIS (defined
below), $\lim_{\Lambda \uparrow V} f_{\Lambda}$ is the limit of the
net $\left\{f_{\Lambda}, \{\Lambda\}_{\Lambda \in \mathcal{S}_{b}, \;
\Lambda \subset V} , \subset \right\}$, for $V$ an infinite interval
of $\mathbb{Z}$. For kernels associated to a specification ,
$\lim_{\Lambda \uparrow V} \gamma_{\Lambda}$ is the limit of the net
$\left\{\gamma_{\Lambda}, \{\Lambda\}_{\Lambda \in \mathcal{S}, \;
\Lambda \subset V} , \subset \right\}$, for $V$ an infinite subset of
$\mathbb{Z}$.  To lighten up formulas involving probability kernels,
we will freely use $\rho(h)$ instead of $E_\rho(h)$ for $\rho$ a
measure on $\Omega$ and $h$ a $\mathcal{F}$-measurable function.  Also
$\rho(\sigma_\Lambda)$ will mean
$\rho(\{\omega\in\Omega:\omega_\Lambda=\sigma_\Lambda\})$ for
$\Lambda\subset\mathbb{Z}$ and $\sigma_\Lambda\in\Omega_\Lambda$.
\bigskip

We start by briefly reviewing the well known notion of specification.
\begin{defi}
A \embf{specification} $\gamma$ on $\left(\Omega,\mathcal{F}\right)$
is a family of probability kernels $\left\{ \gamma_{\Lambda}
\right\}_{\Lambda \in \mathcal{S}}$, $\gamma_{\Lambda} : \mathcal{F}
\times \Omega \rightarrow [0,1]$ such that for all $\Lambda$ in
$\mathcal{S},$
\begin{fleuvea}
\item For each $\displaystyle{A \in \mathcal{F}, \; \gamma_{\Lambda}(A
\mid \cdot ) \in \mathcal{F}_{\Lambda^{c}}.}$
\item For each $\displaystyle{B \in \mathcal{F}_{ \Lambda^{c}} \text{
and } \omega \in \Omega, \; \gamma_{\Lambda}(B \mid \omega) =
\ind{B}(\omega).}$
\item For each $\displaystyle{\Delta \in \mathcal{S}: \Delta \supset
\Lambda,}$ \equa{\label{gibbs5} \gamma_{\Delta}\,\gamma_{\Lambda} \;=\;
\gamma_{\Delta}.}
\end{fleuvea}
The specification is:
\begin{fleuvei}
\item \embf{Continuous} on $\Omega$ if for all $\Lambda
\in \mathcal{S}$ and all $\sigma_\Lambda\in\Omega_\Lambda$ the
functions
\equa{ \label{r.cont}
\Omega\ni\omega \;\longrightarrow\; \gamma_{\Lambda} \left(
\sigma_{\Lambda} \mid \omega \right)
}
are continuous.
\item \embf{Non-null} on $\Omega$ if
$\gamma_{\Lambda}(\omega_{\Lambda} \mid \omega)>0$ for each $\omega
\in \Omega$ and $\Lambda \in \mathcal{S}.$
\end{fleuvei}
\end{defi}
Property c) is usually referred to as \emph{consistency}. There and in the sequel we adopt the standard notation for
composition of probability kernels (or of a probability kernel with a measure).  For instance, \reff{gibbs5} is
equivalently to
$$ \iint h(\xi)\, \gamma_{\Lambda}(d \xi \mid \sigma)
\,\gamma_{\Delta}(d \sigma \mid \omega) \;=\;
\int h(\sigma) \,\gamma_{\Delta} (d \sigma \mid \omega)$$
for each $\mathcal{F}$-measurable function $h$ and configuration
$\omega \in \Omega$.

\begin{remas1}{2}
\item A \emph{Markov specification of range $k$} corresponds to the
particular case in which the applications \reff{r.cont} are in fact
$\mathcal{F}_{\partial_k\Lambda}$-measurable, where
$\partial_k\Lambda=\{i\in \Lambda^{\rm c} : |i-j|\le k \hbox{ for some
} j\in\Lambda\}$.
\item In the sequel, we find useful to consider also the natural
extension of the kernels $\gamma_{\Lambda}$ to functions $\mathcal{F} \times \mathcal{A}^{\mathbb{Z}} \rightarrow
[0,1]$ such that $\gamma_{\Lambda}(\, \cdot \mid \omega)=0$ if $\omega \notin \Omega$. We shall not distinguish
notationally both types of kernels.
\end{remas1}
\begin{defi}\label{d.con}
A probability measure $\mu$ on $(\Omega, \mathcal{F})$ is said to be
\embf{consistent} with a specification $\gamma$ if
\equa{\label{gibbs9} \mu \,\gamma_{\Lambda}\;=\;\mu \quad \forall \;
\Lambda \in \mathcal{S}.}
The family of these measures will be denoted $\mathcal{G}(\gamma)$.
\end{defi}
\begin{remas1}{2}
\item A \emph{Markov field of range $k$} is a measure consistent with
a Markov specification of range $k$.
\item A \embf{Gibbs measure} on $(\Omega, \mathcal{F})$ is a measure
$\mu$ consistent with a specification that is continuous and non-null on $\Omega$.  The SRB measures (Bowen,
1975)\nocite{bow75} are particular one-dimensional examples.
\end{remas1}
\bigskip

We now introduce the analogous notion for processes.  Due to the
nature of the defining transition probabilities, the corresponding
finite-region kernels must apply to functions measurable only with
respect to the region and its past.  Furthermore, finite intervals
already suffice.

\begin{defi}\label{lis1}
A \embf{left interval-specification (LIS)} $f$ on $(\Omega, \mathcal{F})$ is a family of probability kernels
$\left\{f_{\Lambda} \right\}_{ \Lambda \in \mathcal{S}_{b}}$, $f_{\Lambda} : \mathcal{F}_{\leq m_\Lambda} \times \Omega
\longrightarrow [0,1]$ such that for all $\Lambda$ in $\mathcal{S}_{b}$,
\begin{fleuvea}
\item For each $\displaystyle{A \in \mathcal{F}_{\leq m_\Lambda},
\; f_{\Lambda}(A \mid \cdot \, )}$ is
$\mathcal{F}_{\Lambda_{-}}$-measurable.
\item For each $\displaystyle{B \in \mathcal{F}_{\Lambda_{-}} \text{
and } \omega \in \Omega, \; f_{\Lambda}(B \mid \omega) =
\ind{B}(\omega).}$
\item For each $\displaystyle{\Delta \in \mathcal{S}_{b} : \Delta
\supset \Lambda,}$
\equa{\label{lis3} f_{\Delta} \,f_{\Lambda} \;=\; f_{\Delta}
\quad \text{over } \mathcal{F}_{\leq m_\Lambda},}
that is, $(f_{\Delta} f_{\Lambda})(h\mid\omega) =
f_{\Delta}(h\mid\omega)$ for each $\mathcal{F}_{\leq \max
\Lambda}$-measurable function $h$ and configuration $\omega \in
\Omega.$

\end{fleuvea}
The LIS is:
\begin{fleuvei}
\item \embf{Continuous} on $\Omega$ if for all $\Lambda \in \mathcal{S}_{b}$
 and all $\sigma_\Lambda\in\Omega_\Lambda$ the functions
\equa{\label{r.conl}
\Omega\ni\omega \;\longrightarrow\; f_{\Lambda} ( \sigma_{\Lambda}
 \mid \omega )
} are continuous.
\item \embf{Non-null} on $\Omega$ if
$ f_{\Lambda}(\omega_{\Lambda} \mid \omega_{\Lambda_{-}})>0$
for all $\Lambda \in \mathcal{S}_{b}$ and $\omega \in
\Omega_{-\infty}^{m_\Lambda}$.
%
%\item \embf{Strongly non-null} on $\Omega$ if for all $\omega \in
%\Omega, \; \Lambda \in \mathcal{S}$, there exists
%$c\left(\Lambda,\omega_{\Lambda}\right)>0$ such that for all $n \geq
%m_{\Lambda}$
%
%\equa{\label{gibbs22}
%\frac{f_{[l_{\Lambda,n}]}\left(\omega_{[l_{\Lambda,n}]} \mid
%\omega_{\Lambda_{-}}
%\right)}{f_{[l_{\Lambda,n}]}\left(\omega_{\Lambda^{c}\cap[l_{\Lambda,n}]}
%\mid \omega_{\Lambda_{-}} \right)} \;\geq\; c\left( \Lambda,
%\omega_{\Lambda}\right)\;.}
%
\item \embf{Weakly non-null} on $\Omega$ if for all $\Lambda \in
\mathcal{S}_b$, there exists a $\sigma_\Lambda \in \Omega_\Lambda$
such that
$f_\Lambda(\sigma_\Lambda\mid\omega_{-\infty}^{l_\Lambda-1})>0$ for
all $\omega_{-\infty}^{l_\Lambda-1}\in\Omega_{-\infty}^{l_\Lambda-1}$
such that $\sigma_\Lambda\,\omega_{-\infty}^{l_\Lambda-1}\in
\Omega_{-\infty}^{m_\Lambda}$.
\end{fleuvei}
\end{defi}
\begin{remas1}{2}
\item A \emph{Markov LIS of range $k$} is a LIS such that each of the
functions \reff{r.conl} is measurable with respect to
$\mathcal{F}_{[l_\Lambda-k,\, l_\Lambda-1]}$.
\item As for specifications, in the sequel we shall not distinguish
notationally the kernels $f_{\Lambda}$ from their extensions on
$\mathcal{F}_{m_\Lambda} \times \mathcal{A}^{\mathbb{Z}}
\rightarrow [0,1]$ such that $f_{\Lambda}(\, \cdot \mid \omega)=0$ if
$\omega \notin \Omega$.
\end{remas1}
\begin{defi}
A probability measure $\mu$ on $(\Omega,\; \mathcal{F} )$ is said to
be \embf{consistent} with a LIS $f$ if for each $\Lambda \in
\mathcal{S}_{b}$
\equa{\mu\, f_{\Lambda} \;= \;\mu \quad\text{over }
\mathcal{F}_{\leq m_\Lambda}.}
The family of these measures will be denoted $\mathcal{G}(f)$.
\end{defi}
\begin{remas1}{2}
\item
A \emph{Markov chain} of range $k$ is a measure consistent with a Markov LIS of range $k$.
\item Measures consistent with general, non-necessarily Markovian LIS
were initially called \emph{Chains with complete connections} by Onicescu and Mihoc (1935)\nocite{onimih35}. These
objects have been reintroduced several times in the literature under a variety of names: \emph{chains of infinite
order} (Harris, 1955\nocite{har55}), \emph{$g$-measures} (Keane, 1972\nocite{kea72}), \emph{uniform martingales}
(=random Markov processes) (Kalikow, 1990\nocite{kal90}), \ldots.
\end{remas1}

Finally, we introduce a strong notion of uniqueness needed in the sequel.
\begin{defi}
\begin{fleuve1}
\item
  A specification $\gamma$ satisfies a \embf{hereditary uniqueness
  condition} (HUC) \embf{for a family} $\mathcal{H}$ of subsets of
  $\mathbb{Z}$ if for all (possibly infinite) sets $V \in \mathcal{H}$
  and all configurations $\omega \in \Omega$, the specification
  $\gamma^{(V,\omega)}$ defined by
\equa{\label{lisspe0}
\gamma_{\Lambda}^{(V,\omega)}(\, \cdot \mid \xi)=\gamma_{\Lambda}( \,
\cdot \mid \omega_{V^{c}} \, \xi_{V})\;,
\quad\forall \, \Lambda \in \mathcal{S},\, \Lambda \subset V, \;
\; \forall\,\omega_{V^{c}} \, \xi_{V} \in \Omega\;,}
admits a unique Gibbs measure. The specification satisfies a HUC if it
satisfies a HUC for $\mathcal{H} = \mathcal{P}(\mathbb{Z})$.
\item
  A LIS $f$ satisfies a \embf{hereditary uniqueness condition} (HUC)
  if for all intervals of the form $V=[i,+\infty[$, $i\in
  \mathbb{Z}$, or $V=\mathbb{Z}$, and all configurations $\omega \in
  \Omega$, the LIS $f^{(V,\omega)}$ defined by
\equa{\label{lis93}f_{\Lambda}^{(V,\omega)}( \, \cdot \mid
\xi)=f_{\Lambda}(\, \cdot \mid \omega_{V_{-}} \, \xi_{V})\;, \quad
\forall \, \Lambda \in \mathcal{S}_{b},\,  \Lambda \subset V,
\; \forall\ \omega_{V^{c}} \, \xi_{V} \in \Omega\;,}
admits a unique consistent chain.
\end{fleuve1}
\end{defi}

\section{Preliminary results}\label{sect2}
Let us summarize a number of useful properties of LIS and specifications. First we introduce functions associated to
LIS singletons.  For a LIS $f$ and a configuration $\omega \in \Omega$, let \equa{ \label{r.func} f_{i}(\omega)
\;\triangleq\; f_{\{i\}}(\omega_{i} \mid \omega_{-\infty}^{i-1}). } In the shift-invariant case, the function $f_0$ is
a $g$-function in the sense of Keane (1972\nocite{kea72}).

The following theorem expresses the equivalence between the
description in terms of LIS and the usual description in terms of
transition probabilities (=LIS singletons).

\begin{theo}[singleton consistency for chains]\label{listh1}
Let $(g_{i})_{i \in \mathbb{Z}}$ be a family of measurable functions
over $(\Omega, \mathcal{F})$ which enjoy the following properties
\begin{fleuvea}
\item Measurability: for every $i$ in $\mathbb{Z}, \; g_{i} \text{ is
    } \mathcal{F}_{\leq i}\text{-measurable}.$
\item Normalization: for every $i$ in $\mathbb{Z}$ and $\omega \in
  \Omega_{-\infty}^{i-1}$,
\equa{
 \sum_{\sigma_{i} \in \mathcal{A}\,:\; \omega \,
\sigma_{i} \in \Omega_{-\infty}^{i}}
  g_{i}\left(\omega \, \sigma_{i}\right)\;=\;1.
}
\end{fleuvea}
Then there exists a unique left interval-specification $f \triangleq
(f_{\Lambda})_{ \Lambda \in \mathcal{S}_{b}}$ such that $ f_{i}=
g_{i}, \text{ for all } i \text{ in } \mathbb{Z}$.  Furthermore:
\begin{fleuvei}
\item $f$ satisfies (in fact, it is defined by) the property
\equa{\label{lis5} f_{[l,m]}=f_{[l,n]}f_{[n+1,m]} \quad \text{ over }
  \mathcal{F}_{\leq n}}
for each $l,m,n \in \mathbb{Z} : l \leq n < m$.
\item $f$ is non-null on $\Omega$ if, and only if, so are the
functions $g_i$, that is, if and only if $g_i(\omega)>0$ for each
$i\in\mathbb{Z}$ and each $\omega\in\Omega_{-\infty}^i$.
\item $f$ is weakly non-null on $\Omega$ if, and only if, so are the
functions $g_i$, that is, if and only if for each $i\in\mathbb{Z}$
there exists $\sigma_i \in \Omega_{\{i\}}$ such that
$g_i(\sigma_i\,\omega_{-\infty}^{i-1})>0$
for all
$\omega_{-\infty}^{i-1}\in\Omega_{-\infty}^{i-1}$
such that
$\sigma_i\,\omega_{-\infty}^{i-1}\in \Omega_{-\infty}^{i}$.
\item $\mathcal{G}(f) = \{ \mu : \mu f_{i}= \mu, \text{ for
    all } i \text{ in } \mathbb{Z} \}$.
\end{fleuvei}
\end{theo}
The proof of this result is rather simple (it is spelled up in Fern\'andez and Maillard, 2003\nocite{fermai03b}).  The
following theorem is the analogous result for specifications.  Let us consider the following functions associated to an
specification $\gamma$: \equa{ \label{r.funs} \gamma_{\Lambda}(\omega) \;\triangleq\; \gamma_{\{\Lambda\}}
(\omega_{\Lambda} \mid \omega_{\Lambda^{\rm c}})\quad,\quad \gamma_i(\omega) \;\triangleq\; \gamma_{\{i\}}(\omega)\;. }

\begin{theo}[singleton consistency for Gibbs measure]\label{speth1}
Let $\left( \rho_{i} \right)_{i \in \mathbb{Z}}$ be a family of
measurable functions over $(\mathcal{A}^{\mathbb{Z}}, \mathcal{F})$
which enjoys the following properties
\begin{fleuvea}
\item Non-nullness on $\Omega$: for every $i$ in $\mathbb{Z}, \;
\rho_{i}(\omega)=0 \Longleftrightarrow \omega \notin \Omega.$
\item Order-consistency on $\Omega$: for every $i,j \in \mathbb{Z}$
and $\omega \in \Omega$,
\equa{\label{r.ord}
\frac{\rho_{i}(\omega)}{\displaystyle{\sum_{\sigma_{i} \in \mathcal{A}}
\rho_{i}\left(\sigma_{i}\,\omega_{\{i\}^{c}} \right) \,
\rho_{j}^{-1}\left(\sigma_{i}\,\omega_{\{i\}^{c}}\right)}} \;=\;
\frac{\rho_{j}(\omega)}{\displaystyle{ \sum_{\sigma_{j} \in
\mathcal{A}} \rho_{j} \left(\sigma_{j}\,\omega_{\{j\}^{c}}\right)
\, \rho_{i}^{-1} \left(\sigma_{j}\,\omega_{\{j\}^{c}}\right)}}.
}
\item Normalization on $\Omega$: for every $i \in \mathbb{Z}$ and
$\omega \in \Omega$,
\equa{
\sum_{\sigma_{i} \in \mathcal{A}} \rho_{i}\left(\sigma_{i}\,
\omega_{\{i\}^{c}} \right)=1.
}
\end{fleuvea}
Then there exists a unique specification $\gamma$ on
$(\Omega, \mathcal{F})$ such that $ \gamma_{i}(\omega)=
\rho_{i}(\omega), \text{ for all } i \text{ in } \mathbb{Z}$.
Furthermore,
\begin{fleuvei}
\item $\gamma$ is non-null on $\Omega$: For each $\Lambda \in
\mathcal{S}$, $\gamma_{\Lambda}(\omega)=0 \Longleftrightarrow \omega
\notin \Omega.$
\item $\gamma$ satisfies an order-independent prescription:
For each $\Lambda, \Gamma \in \mathcal{S}$ with $\Gamma \subset
\Lambda^{\rm c}$
\equa{
\label{gibbs37}
\gamma_{\Lambda\cup\Gamma}(\omega) \;=\;
\frac{\gamma_{\Lambda}(\omega)}{\displaystyle{\sum_{\sigma_{\Lambda}}
\gamma_{\Lambda} \left( \sigma_{\Lambda}\,\omega_{\Lambda^{\rm c}}
 \right) \, \gamma_{\Gamma}^{-1}\left(
\sigma_\Lambda\,\omega_{\Lambda^{\rm c}}\right)}}
}
for all $\omega \in \Omega$.
\item $\mathcal{G}(\gamma) = \left\{ \mu : \mu\, \gamma_{i} = \mu \text{ for
all } i \in \mathbb{Z} \right\}$.
\end{fleuvei}
\end{theo}
This result will be proved in Appendix \ref{A2} in a more general
setting.
\bigskip

For completeness, we list now several, mostly well known, sufficient conditions for hereditary uniqueness.  They refer
to different ways to bound continuity rates of transition kernels. We start with the relevant definitions.

\begin{defi}
\begin{fleuvei}
\item The \embf{$k$-variation} of a $\mathcal{F}_{\{i\}}$-measurable
function $f_i$ is defined by
$$ \var{k}(f_i) \;\triangleq\; \sup
\Bigl\{\left| f_i(\omega_{-\infty}^i) - f_i(\sigma_{-\infty}^i)\right|
: \omega_{-\infty}^i, \sigma_{-\infty}^i \in \Omega_{-\infty}^i \,,
\, \omega_{i-k}^i = \sigma_{i-k}^i\Bigr\}\;.$$
\item The \embf{interdependence coefficients} for a
family of probability kernels $\pi=\left( \pi_{\{i\}}\right)_{i \in
\mathbb{Z}}, \;\pi_{\{i\}} : \mathcal{F}_{i} \times \Omega \rightarrow
[0,1]$ are defined by
\equa{\label{int7} C_{ij}(\pi) \;\triangleq\;
\sup_{\substack{\xi, \eta \in \Omega\\ \xi \stackrel{\neq j}{=} \eta}}
\left\| \stackrel{\circ}{\pi}_{\{i\}} (\, \cdot \mid \xi) -
\stackrel{\circ}{\pi}_{\{i\}}(\, \cdot \mid \eta )\right\|}
for all $i,j \in \mathbb{Z}$. Here we use the variation norm and
$\stackrel{\circ}{\pi}_{\{i\}}$ is the projection of $\pi_{\{i\}}$
over $\{i\}$ that is $\stackrel{\circ}{\pi}_{\{i\}} \left( A \mid
\omega\right) \triangleq \pi_{\{i\}}\left( \left\{ \sigma_{i} \in A
\right\} \mid \omega \right)$ for all $A \in \mathcal{F}_{i}$ and
$\omega \in \Omega$.
\end{fleuvei}
\end{defi}

A LIS $f$ on $(\Omega, \mathcal{F})$ satisfies a HUC if it satisfies one of the following conditions:
\begin{itemize}
\item \embf{Harris} (Harris 1955; Coelho and Quas, 1998)\nocite{har55}\nocite{coequa98}: The LIS $f$ is stationary,
weakly non-null on $\Omega$ and
$$\sum_{n \geq 1} \prod_{k=1}^{n} \left( 1 - \frac{|\mathcal{A}|}{2} \var{k}(f_{0}) \right) = + \infty.$$
\item \embf{Berbee}\nocite{ber87} (1987): The LIS $f$ is stationary, non-null and $$ \sum_{n \geq 1} \exp \left( -
\sum_{k=1}^{n} \var{k}(\log f_0) \right) = + \infty.$$
\item \embf{Stenflo}\nocite{ste03} (2002): The LIS $f$ is stationary, non-null and $$\sum_{n \geq 1} \prod_{k=1}^{n}
\Delta_{k}(f_0) = + \infty$$ where $\Delta_{k}(f_0) \triangleq \inf \{ \sum_{\omega_{0} \in \mathcal{A}} \min \left(f_0
\left(\omega_{-\infty}^{0}\right),f_0 \left(\sigma_{-\infty}^{-1} \, \omega_{0}\right) \right) :
\omega_{-k}^{-1}=\sigma_{-k}^{-1} \}$.
\item \embf{Johansson and {\"O}berg}\nocite{johobe02} (2002): The LIS $f$ is stationary, non-null and $$\sum_{k \geq 0}
\var{k}^{2}(\log f_0)< + \infty.$$
\item \embf{One-sided Dobrushin} (Fern\`andez and Maillard, 2003): For each $\displaystyle{ i \in \mathbb{Z}, \;
\sum_{j<i} \, C_{ij}(f) < 1}$ and $f$ is continuous on $\Omega$.
\item \embf{One-sided boundary-uniformity} (Fern\`andez and Maillard, 2003): There exists a constant $K>0$ so that for
every cylinder set $A=\{x_{l}^{m}\} \in \Omega_{l}^{m}$ there exists an integer $n$ such that
\equa{\label{lis91} f_{[n,m]}(A \mid \xi)\geq K \, f_{[n,m]}(A \mid \eta) \quad \text{for all } \xi, \eta \in
\Omega\;.}
\end{itemize}

The last two conditions, proven in a companion paper (Fern\'andez and Maillard, 2003)\nocite{fermai03b}, are in fact
adaptations of the following well known criteria for specifications.

A specification $\gamma$ on $(\Omega, \mathcal{F})$ satisfies a HUC if
it satisfies one of the following conditions:
\begin{itemize}
\item \embf{Dobrushin} (1968), Lanford
(1973)\nocite{dob68}\nocite{lan73}: $\displaystyle{ \sup_{i \in
\mathbb{Z}} \sum_{j \in \mathbb{Z}} \, C_{ij}( \gamma) < 1}$ and
$\gamma$ continuous.
\item \embf{Georgii} (1974)\nocite{geo74} \embf{boundary-uniformity}: There
exists a constant $K>0$ so that for every cylinder set $A \in
\mathcal{F}$ there exists $\Lambda \in \mathcal{S}_{b}$ such that
$$\gamma_{\Lambda}(A \mid \xi)\geq K \, \gamma_{\Lambda}(A \mid \eta)
\quad \text{for all } \xi, \eta \in \Omega\;.$$
\end{itemize}
We remark that the conditions involve no non-nullness assumption.

\section{Main results}
For a LIS $f$ on $\Omega$ let us denote, for each $\Lambda \in
\mathcal{S}, \; k \geq m_{\Lambda}$ and $\omega \in
\Omega_{-\infty}^{k}$,
\equa{\label{def5} c_{\Lambda}^{\omega}(f_{k}) \;\triangleq \;\inf_{\sigma_{\Lambda} \in \Omega_{\Lambda} }
\left\{f_{k}\left(\sigma_{\Lambda}\, \omega_{]-\infty,k] \setminus \Lambda} \right) : \sigma_{\Lambda} \,
\omega_{]-\infty,k] \setminus \Lambda} \in \Omega_{-\infty}^{k} \right\}}
and
\equa{\label{def7}\delta_{\Lambda}^{\omega}(f_{k}) \;\triangleq\; \sum_{j
\in \Lambda} \sup \left\{ \left| f_{k} \left(
\omega_{-\infty}^{k}\right) - f_{k} \left( \sigma_{-\infty}^{k}
\right) \right| : \sigma \in \Omega_{-\infty}^{k}, \; \sigma
\stackrel{\neq j}{=} \omega \right\}.}
Similarly, for a specification $\gamma$ on
$\Omega$, let us denote, for each $\omega \in \Omega$ and each $k,j \in
\mathbb{Z}$
\equa{\label{def9} c_{j}^{\omega}(\gamma_{k})\; \triangleq\;
\min_{\sigma_{j} \in \Omega_{\{j\}}} \left\{
\gamma_{k}\left(\sigma_{j}\,\omega_{\{j\}^{c}} \right) :
\sigma_{j}\,\omega_{\{j\}^{c}} \in \Omega \right\}}
and
\equa{\label{def13}\delta_{j}^{\omega}\left(\gamma_{k} \right)
\;\triangleq\; \sup \left\{ \left| \gamma_{k}(\omega) - \gamma_{k}(
\sigma) \right| : \sigma \in \Omega, \; \sigma \stackrel{\neq j}{=}
\omega \right\}\;.}

\begin{defi}
\begin{fleuvei}
\item A LIS $f$ on $\Omega$ is said to have a \embf{good future (GF)}
if it is non-null on $\Omega$ and for each $\Lambda \in
\mathcal{S}$, there exists a sequence $\left\{
\varepsilon_{k}^{\Lambda} \right\}_{k \in \mathbb{N}}$ of positive
numbers such that $\sum_{k} \varepsilon_{k}^{\Lambda} < +\infty$ for
which
\equa{\label{lisspe1}
\sup_{\omega \in \Omega_{-\infty}^k}\,
c_{\Lambda}^{\omega}(f_{k})^{-1} \,
\delta_{\Lambda}^{\omega}(f_{k}) \;\leq\; \varepsilon_{k}^{\Lambda}}
for each $k \geq m_{\Lambda}$.
\item A LIS $f$ on $\Omega$ is said to have an
\embf{exponentially-good future (EGF)} if it is non-null on
$\Omega$ and there exists a real $a>1$ such that
\equa{
\limsup_{k\to\infty} \,a^{|k-j|}\, \sup_{\omega \in
\Omega_{-\infty}^k}\, c_{j}^{\omega}(f_{k})^{-1} \,
\delta_{j}^{\omega}(f_{k}) \;<\;\infty
}
for all $j \in \mathbb{Z}$.
\item A specification $\gamma$ on $\Omega$ is said to have an
\embf{exponentially-good future (EGF)} if it is non-null on $\Omega$
and there exists a real $a>1$ such that
\equa{
\limsup_{k\to\infty} \,a^{|k-j|}\; \sup_{\omega \in \Omega}\;
c_{j}^{\omega}(\gamma_{k})^{-1} \, \delta_{j}^{\omega}
\left(\gamma_{k}\right) \;<\; \infty
}
for all $j \in \mathbb{Z}$.
\end{fleuvei}
\end{defi}

\begin{defi}
Let us introduce the following sets.
\medskip

\hspace{1cm}
\begin{minipage}{7cm}
$\Theta \triangleq \left\{\hbox{LIS }f \hbox{continuous and
non-null on } \Omega\right\}$\\
$\Pi \triangleq \left\{\hbox{specifications }\gamma
\hbox{ continuous and non-null on } \Omega\right\}$,
$\Theta_{1} \triangleq \left\{ f \in \Theta : f \text{ has a GF}
\right\}$,\\
$\Pi_{1} \triangleq \left\{\gamma \in \Pi : \left| \mathcal{G}
(\gamma) \right| = 1 \right\}$,\\
$\Pi_{2} \triangleq \left\{ \gamma \in \Pi : \gamma \text{ satisfies a
HUC over all } [i,+\infty[, \; i \in \mathbb{Z} \right\}$,\\
$\Theta_{2} \triangleq \left\{ f \in \Theta : f \text{ satisfies a
HUC}\right\}$,\\
$\Theta_{3} \triangleq \left\{ f \in \Theta : f \text{ has an EGF}
\right\}$,\\
$\Pi_{3} \triangleq \left\{ \gamma \in \Pi : \gamma \text{ has an EGF}
\right\}$.
\end{minipage}
\end{defi}
%\bigskip

We remark that each of the LIS or specifications of any of the
preceding sets has at least one consistent measure.  This is because
the (interesting part of) the configuration space is compact and the
LIS or specifications are assumed to be continuous.  Indeed, as the
space of probability measures on a compact space is weakly compact,
every sequence of measures
$\gamma_{\Lambda_n}(\,\cdot\,\mid\omega^{\{n\}})$ or
$f_{\Lambda_n}(\,\cdot\,\mid\omega^{\{n\}})$, for $(\Lambda_n)$ an
exhausting sequence of regions and $(\omega^{\{n\}})$ a sequence of
configurations, has a weakly convergent subsequence.  By continuity of
the transitions the limit is respectively a Gibbs measure or a
consistent chain.
\medskip

Consider the function
\equa{\label{comp2}
F_{\Lambda,n}\left(\omega_{\Lambda} \mid
\omega\right) \;\triangleq\; \frac{\displaystyle{f_{[l_{\Lambda},n]}
\left( \omega_{l_{\Lambda}}^{n} \mid \omega_{\Lambda_{-}} \right)}}{
\displaystyle{ f_{[l_{\Lambda},n]} \left( \omega_{\Lambda^{c} \cap
[l_{\Lambda},n]} \mid \omega_{\Lambda_{-}} \right)}}}
for all $\Lambda \in \mathcal{S}$, $n \geq m_{\Lambda}$ and $\omega
\in \Omega$.  The continuity of $f$ implies that the functions
$F_{\Lambda,n}(\omega_\Lambda \mid \cdot \,)$ are continuous on
$\Omega_{\Lambda^{c}}$ for each $\omega_{\Lambda} \in
\Omega_{\Lambda}$. We use these functions to introduce the map
$$b : \Theta_{1} \rightarrow \Pi\;, \; f \mapsto \gamma^{f}$$
defined by
\equa{\label{comp3} \gamma_{\Lambda}^{f}\left( \omega_{\Lambda} \mid
\omega \right)\;\triangleq\; \lim_{n \rightarrow +\infty}
F_{\Lambda,n}(\omega_{\Lambda} \mid \omega)}
for all $\Lambda \in \mathcal{S}$ and $\omega \in \Omega$.

\begin{theo}[LIS $\rightsquigarrow$ specification]\label{th1}
\begin{fleuve1}
\item The map $b$ is well defined. That is, for $f \in \Theta_{1}$
\begin{fleuvea}
\item the limit \reff{comp3} exists for all $\Lambda \in \mathcal{S}$
  and $\omega \in \Omega$.
\item $\gamma^{f}$ is a specification on $(\Omega, \mathcal{F})$.
\item $\gamma^{f} \in \Pi$.
\end{fleuvea}
\item
\begin{fleuvea}
\item For each (finite or infinite) interval $V$ and each $\omega \in
  \Omega, \; \mathcal{G}\left(f^{(V,\omega)} \right) \subset
  \mathcal{G}\left(\left(\gamma^{f} \right)^{(V,\omega)} \right)$.
\item For $f \in b^{-1}(\Pi_{1}), \; \mathcal{G}(f) =
  \mathcal{G}(\gamma^{f}) = \left\{ \mu^{f}\right\}$, where $\mu^{f}$
  is the only chain consistent with $f$.
\item The map $b$ restricted to $b^{-1}\left( \Pi_{1} \right)$ is one-to-one.
\end{fleuvea}
\end{fleuve1}
\end{theo}
\begin{rema}
Since for all $k \geq 1, \, \var{k}(g) \geq \delta_{k}(g)$, $\Theta_{1}$ includes the set of stationary non-null LIS
with summable variation.
%The latter is the setting Ledrappier
%(1974)\nocite{led74}, Walters (1975)\nocite{wal75},Lalley
%(1986)\nocite{lal86}, Berbee (1987)\nocite{ber87}, Bressaud et al
%(1999)\nocite{bfg2}, \ldots.
\end{rema}
Consider now the map
\equa{\label{comp4}
c : \Pi_{2} \rightarrow \Theta_{2}\; , \; \gamma
\mapsto f^{\gamma}}
defined by
\equa{\label{comp5} f_{\Lambda}^{\gamma}\left( A \mid
    \omega_{\Lambda_{-}} \right)\;\triangleq\; \lim_{k \rightarrow
    +\infty} \gamma_{\Lambda \cup \Lambda_{+}^{(k)}}\left( A \mid
    \omega \right)}
for all $\Lambda \in \mathcal{S}_{b}, \; A \in \mathcal{F}_{\Lambda}$,
and $\omega \in \Omega$ for which the limit exists.
\begin{theo}[specification $\rightsquigarrow$ LIS]\label{th2}
\begin{fleuve1}
\item The map $c$ is well defined. That is, for $\gamma \in \Pi_{2}$
\begin{fleuvea}
\item the limit \reff{comp5} exists for all $\Lambda \in
\mathcal{S}_{b}, \; A \in \mathcal{F}_{\Lambda}, \;
\omega_{\Lambda_{-}} \in \Omega_{\Lambda_{-}}$ and is independent of
$\omega_{\Lambda_{+}}$.
\item $f^{\gamma}$ is a LIS on $(\Omega, \mathcal{F})$.
\item $f^{\gamma} \in \Theta_{2}$.
\end{fleuvea}
\item
\begin{fleuvea}
\item $\mathcal{G}(f^{\gamma}) = \mathcal{G}(\gamma) = \left\{
\mu^{\gamma}\right\}$, where $\mu^{\gamma}$ is the only Gibbs measure
consistent with $\gamma$.
\item The map $c$ is one-to-one.
\end{fleuvea}
\end{fleuve1}
\end{theo}
In addition a LIS of the form $f^{\gamma}$ satisfies the following
properties.
\begin{theo}\label{th3}
Let $\gamma \in \Pi_{2}$.
\begin{fleuvea}
\item If $\gamma$ satisfies Dobrushin uniqueness condition, then so
does $f^{\gamma}$.
\item If $\gamma$ satisfies the boundary-uniformity uniqueness
condition, then so does $f^{\gamma}$.
\end{fleuvea}
\end{theo}
\begin{theo}[Continuity rates]\label{th5}
Let $\omega \in \Omega$ and $j \in \mathbb{Z}$.
\begin{fleuve1}
\item For $f \in \Theta_{1}$ and $\Lambda\in\mathcal{S}$
\begin{fleuvea}
\item if $j>m_\Lambda$ then
$\displaystyle{\delta_{j}^{\omega}\left(\gamma_{\Lambda}^{f} \right)
\leq 2 \sum_{i \geq j} c_{k}^{\omega} (f_{i})^{-1} \,
\delta_{k}^{\omega}(f_{i}).}$
\item if $j<l_\Lambda$ then
$\displaystyle{\delta_{j}^{\omega}(\gamma^f_{\Lambda})
\leq 1 - \prod_{i=l_\Lambda}^{+\infty}\frac{ 1 -
c_{j}^{\omega}(f_{i})^{-1} \,
\delta_{j}^{\omega}( f_{i})}{1 + c_{j}^{\omega}(f_{i})^{-1} \,
\delta_{j}^{\omega}( f_{i})}.}$
\end{fleuvea}
\item For $\gamma \in \Pi_{2}$, $\Lambda\in\mathcal{S}_{b}$
 and $j<l_\Lambda$, $\displaystyle{
\delta_{j}^{\omega}(f_{\Lambda}^{\gamma})\leq 1 -
\prod_{i=l_\Lambda}^{+\infty}
\frac{1-c_{j}^{\omega}(\gamma_{k})^{-1} \,
\delta_{j}^{\omega}(\gamma_{i})}{1+c_{j}^{\omega}(\gamma_{k})^{-1} \,
\delta_{j}^{\omega}(\gamma_{i})}.}$
\end{fleuve1}
\end{theo}
Under suitable conditions the maps $b$ and $c$ are reciprocal.
\begin{theo}[LIS $\leftrightsquigarrow$ specification]\label{th4}
\begin{fleuvea}
\item $b \circ c = \Id \quad \text{over } c^{-1}(\Theta_{1})$ and
$\mathcal{G}(f^{\gamma}) = \mathcal{G}(\gamma) = \left\{
\mu^{\gamma}\right\}$.
\item $c \circ b = \Id \quad \text{over } {b^{-1}(\Pi_{2})}$ and
$\mathcal{G}(\gamma^{f}) = \mathcal{G}(f) = \left\{ \mu^{f}\right\}$.
\item $b$ and $c$ establish a one-to-one correspondence between
$\Theta_{3}$ and $\Pi_{3}$ that preserves the consistent measure.
\end{fleuvea}
\end{theo}

We remark that $\Theta_3$ includes the well studied processes with
Hold\"erian transition rates (see, for instance, Lalley,
1986\nocite{lal86}, or Keller, 1998\nocite{kel98}).  Part (c) of the
theorem shows, in particular, the equivalence between such processes
and Bowen's Gibbs measures.
\section{Proofs}
We start with a collection of results used for several proofs.
\begin{lemm} \label{r.lem1}
  Consider $\Lambda \in \mathcal{S}, \; \Delta \subset \Lambda^{c}$
  and $\pi$ a probability kernel over $\left(\mathcal{F}_{\Lambda}
    \otimes \mathcal{F}_{\Delta} \,,\,\Omega_{\Lambda} \times
    \Omega_{\Delta} \right)$ such that $\pi\left(A_{\Delta} \mid \cdot
    \, \right) = \ind{A_{\Delta}}(\, \cdot \,), \; \forall \,
  A_{\Delta} \in \mathcal{F}_{\Delta}$. Then, for all $\omega \in
  \Omega_{\Lambda} \times \Omega_{\Delta}$,
  $$\pi(\, \cdot \mid \omega) \;=\;
  [\stackrel{\circ}{\pi}_\Lambda(\, \cdot \mid \omega) \otimes
  \delta_{\omega_{\Delta}}]( \, \cdot \,)$$
  where
  $\stackrel{\circ}{\pi}_\Lambda(\, \cdot \mid \omega)$ is the
  restriction of $\pi(\, \cdot \mid \omega)$ to
  $\mathcal{F}_{\Lambda}$ and $\delta_{\omega_{\Delta}}$ is the Dirac
  mass at $\omega_{\Delta}$.
\end{lemm}
\proof If $A=A_{\Lambda} \times A_{\Delta} \in \mathcal{F}_{\Lambda} \otimes \mathcal{F}_{\Delta}$ and $\omega \in
\Omega_{\Lambda} \times \Omega_{\Delta}$
\[
 \pi\left (A_\Delta \bigm| \omega \right) \;=\; \pi\left( A_{\Lambda} \times A_{\Delta} \bigm| \omega \right) + \pi\left(
A_{\Lambda}^{c} \times A_{\Delta} \bigm| \omega \right)
\]
and \[
 \pi\left (A_\Lambda \bigm| \omega \right) \;=\; \pi\left( A_{\Lambda} \times A_{\Delta} \bigm| \omega \right) + \pi\left(
A_{\Lambda} \times A_{\Delta}^{c} \bigm| \omega \right) \;.
\]
Hence,
\[
 \pi\left( A_{\Lambda} \times A_{\Delta} \right) \;\leq\;
 \pi\left( A_{\Lambda}\right) \, \wedge \, \pi\left( A_{\Delta} \right) \;\leq\;
\pi \left( A_{\Lambda} \right) \ind{A_{\Delta}}.
\]
Analogously, \[ \pi\left( A_{\Lambda} \times A_{\Delta}^{c} \right) \;\leq\; \pi \left( A_{\Lambda} \right)
\ind{A_{\Delta}^{c}}\;.
\]
On the other hand,
$$\pi\left( A_{\Lambda} \times A_{\Delta} \right)
+ \pi\left( A_{\Lambda} \times A_{\Delta}^{c} \right) \;=\; \pi \left( A_{\Lambda} \right) \ind{A_{\Delta}} + \pi
\left( A_{\Lambda} \right) \ind{A_{\Delta}^{c}}\;.$$ The last three displays imply
$$\pi(A ) = \pi\left(A_{\Lambda}\right)
\; \ind{A_{\Delta}} \quad \Longleftrightarrow \quad
%\text{if and only if} \quad
\pi \left(A_{\Delta} \right) = \ind{A_{\Delta}}\;.\ \square$$

In particular, LIS and specifications are completely defined by the
families of their restrictions.

\begin{prop}\label{l.r.2}
Let $\omega \in \Omega$, $\Lambda \in \mathcal{S}$ and
$n\in\mathbb{Z}$, $n\ge m_\Lambda$.
\begin{fleuvea}
\item For any $\beta \in \Omega$,
\equa{
  \label{eq:2}
f_{n+1} \left( \omega_{n+1} \mid \omega_{]-\infty,n]}
\right) \;{\le\atop\ge}\;
f_{n+1} \left( \omega_{n+1} \mid \beta_{\Lambda} \,
\omega_{]-\infty,n] \setminus \Lambda}\right)  \,
\Bigl[1\,\pm\, {\delta_{\Lambda}^{\omega}(f_{n+1})\over
c_{\Lambda}^{\omega}(f_{n+1})}\Bigr]\;.
}
\item For $j< l_\Lambda$ and $\sigma \in \Omega$ with $\sigma
\stackrel{\neq j}{=} \omega$
\eqna{
\label{eq:rr.1}
\lefteqn{
F_{\Lambda,n}(\omega_{\Lambda} \mid \sigma) -
F_{\Lambda,n}(\omega_{\Lambda} \mid \omega)
\quad{\le\atop\ge}}\nonumber\\[5pt]
& & \pm \ \biggl[ 1 - \prod_{i=l_\Lambda}^{n}\frac{ 1 -
c_{j}^{\omega}(f_{i})^{-1} \,
\delta_{j}^{\omega}( f_{i})}{1 + c_{j}^{\omega}(f_{i})^{-1} \,
\delta_{j}^{\omega}( f_{i})}\biggr]
\,\times\,\left\{{F_{\Lambda,n}(\omega_{\Lambda} \mid \sigma) \atop
F_{\Lambda,n}(\omega_{\Lambda} \mid \omega) }\right.
\;.
}
\item
\equa{
  \label{eq:3}
F_{\Lambda,n+1}(\omega_{\Lambda} \mid \omega)
\;{\le\atop\ge}\;
F_{\Lambda,n}(\omega_{\Lambda} \mid \omega) \,
\Bigl[1\,\pm\, {\delta_{\Lambda}^{\omega}(f_{n+1})\over
c_{\Lambda}^{\omega}(f_{n+1})}\Bigr]
}
\end{fleuvea}
\end{prop}

\proof
If we telescope $ f_{n+1}
\left( \omega_{n+1} \mid \xi_{\Lambda} \, \omega_{]-\infty,n]
    \setminus \Lambda} \right) - f_{n+1} \left( \omega_{n+1} \mid
  \eta_{\Lambda} \, \omega_{]-\infty,n] \setminus \Lambda}\right)$,
transforming $\xi_{\Lambda}$ into $\eta_{\Lambda}$ letter by letter,
we have
\[
\label{lisspe6} f_{n+1} \left( \omega_{n+1} \mid
    \xi_{\Lambda} \, \omega_{]-\infty,n] \setminus \Lambda} \right) -
  f_{n+1} \left( \omega_{n+1} \mid \eta_{\Lambda} \, \omega_{
      ]-\infty,n] \setminus \Lambda}\right) \;\leq\;
  \delta_{\Lambda}^{\omega}(f_{n+1})
\]
for any $\xi_\Lambda,\eta_\Lambda\in\Omega_\Lambda$.  To obtain {\bf \reff{eq:2}} we simply use this inequality twice,
assigning $\omega_\Lambda$ to $\xi_\Lambda$, $\beta_\Lambda$ to $\eta_\Lambda$, and vice versa.
\smallskip

To prove {\bf \reff{eq:rr.1}} we use definition \reff{comp2}
\[
\label{eq:spe5} F_{\Lambda,n}(\omega_{\Lambda} \mid \omega) =
  \frac{f_{[l_{\Lambda},n]} \left(\omega_{l_{\Lambda}}^{n} \mid
      \omega_{\Lambda_{-}} \right)} {
    \displaystyle{\sum_{\beta_{\Lambda}}  f_{[l_{\Lambda},n]}
      \left(\beta_{\Lambda} \, \omega_{\Lambda^{c} \cap
          [l_{\Lambda},n]} \mid \omega_{\Lambda_{-}}\right) }}
\]
and the factorization
\[
\label{exp2.1}
f_{[k,n]}\left( \omega_{k}^{n} \mid \omega_{-\infty}^{k-1} \right)=
\prod_{i=k}^{n} f_{i}\left( \omega_{i} \mid
\omega_{-\infty}^{i-1}\right)\;.
\]
We then apply inequalities \reff{eq:2} to bound each of the factors by
similar factors with conditioning configuration $\sigma$.
\smallskip

To obtain {\bf \reff{eq:3}} we apply the LIS-reconstruction
formula \reff{lis5} with $m = n+1$ which yields
\[
\label{lisspe5} F_{\Lambda,n+1}(\omega_{\Lambda} \mid \omega) =
  \frac{f_{[l_{\Lambda},n]} \left(\omega_{l_{\Lambda}}^{n} \mid
      \omega_{\Lambda_{-}} \right) \, f_{n+1} \left( \omega_{n+1} \mid
      \omega_{]-\infty,n]} \right)}{
    \displaystyle{\sum_{\beta_{\Lambda}} \bigg[ f_{[l_{\Lambda},n]}
      \left(\beta_{\Lambda} \, \omega_{\Lambda^{c} \cap
          [l_{\Lambda},n]} \mid \omega_{\Lambda_{-}}\right) \, f_{n+1}
      \left( \omega_{n+1} \mid \beta_{\Lambda} \, \omega_{]-\infty,n]
          \setminus \Lambda} \right)\bigg]}}\;.
\]
In the denominator, only $\beta_{\Lambda}$ with $f_{n+1} \left(
  \omega_{n+1} \mid \omega_{]-\infty,n] \setminus \Lambda} \,
  \beta_{\Lambda}\right)\neq 0$ contribute. We use inequalities
  \reff{eq:2} for these.  $\square$

\begin{lemm}\label{lemproof}
  Let $f \in \Theta_{1}, \; \omega \in \Omega$ and $\Lambda \in
  \mathcal{S}$. Consider the sequence
  $F_{\Lambda,n}\left(\omega_{\Lambda} \mid \omega\right), \, n \geq
  m_{\Lambda}$ defined by \reff{comp2}. Then the limit \reff{comp3}
  exits and satisfies \equa{\label{lisspe51} \left|
      \gamma_{\Lambda}^{f} \left( \omega_{\Lambda} \mid
        \omega_{\Lambda^{c}}\right) -
      F_{\Lambda,n}\left(\omega_{\Lambda} \mid \omega\right) \right|
    \leq \sum_{k \geq n+1} c_{\Lambda}^{\omega}(f_{k})^{-1} \,
    \delta_{\Lambda}^{\omega}(f_{k}).}
\end{lemm}
\proof
From \reff{eq:3}, plus the fact that $F_{\Lambda,n}
\in [0,1] \; \forall \, n \geq m,$ we obtain
$$\left| F_{\Lambda,n+1}\left(\omega_{\Lambda} \mid \omega \right) -
  F_{\Lambda,n}\left(\omega_{\Lambda} \mid \omega \right) \right| \leq
c_{\Lambda}^{\omega}(f_{n+1})^{-1} \,
\delta_{\Lambda}^{\omega}(f_{n+1}).$$
Therefore the summability of
$c_{\Lambda}^{\omega}(f_{k})^{-1} \, \delta_{\Lambda}^{\omega}(f_{k})$
for $k \geq m$ implies the summability of the sequence $\left|
  F_{\Lambda,k+1}\left(\omega_{\Lambda} \mid \omega \right) -
  F_{\Lambda,k}\left(\omega_{\Lambda} \mid \omega \right) \right|$. In
particular $\big( F_{\Lambda,n}(\omega_{\Lambda} \mid \omega )
\big)_{n \geq m}$ is a Cauchy sequence so the limit
  $\lim_{n \rightarrow +\infty} F_{\Lambda,n}\left(\omega_{\Lambda} \mid
  \omega \right) \triangleq \gamma_{\Lambda}^{f} \left(
  \omega_{\Lambda} \mid \omega \right)$ exits and satisfies
\reff{lisspe51} for each $\omega \in \Omega. \quad \Box$
\subsection{LIS $\rightsquigarrow$ specification}

\paragraph{Proof of Theorem \protect \ref{th1}}\mbox{}
\smallskip

Lemma \reff{lemproof} proves Item {\bf 1) (a)}.
\smallskip

To prove item {\bf 1) (b)}, we observe that $\gamma_{\Lambda}^{f}(A
\mid \cdot \, )$ is clearly $\mathcal{F}_{\Lambda^{c}}$-measurable for
every $\Lambda \in \mathcal{S}$ and every $A \in \mathcal{F}$.
Moreover condition (b) of Definition \ref{lis1} together with the
presence of the indicator function $\ind{\omega_{\Lambda^{c} \cap
    [l_{\Lambda},m_{\Lambda}]}}$ in the denominator of \reff{comp2}
imply that $\gamma_{\Lambda}^{f}(B \mid \cdot \, )= \ind{B}(\, \cdot
\,)$ for every $\Lambda \in \mathcal{S}$ and every $B \in
\mathcal{F}_{\Lambda^{c}}$. Therefore it suffices to show that
\equa{\label{lisspe53} \sum_{\omega_{\Delta \setminus
\Lambda}}\gamma_{\Lambda}^{f}\left( \omega_{\Lambda} \mid
\omega_{\Lambda^{c}} \right) \gamma_{\Delta}^{f}\left( \omega_{\Delta
\setminus \Lambda} \mid \omega_{\Delta^{c}} \right) \;=\;
\gamma_{\Delta}^{f}\left( \omega_{\Lambda} \mid
\omega_{\Delta^{c}}\right)}
for each $\Lambda, \Delta \in \mathcal{S}$ such that $\Lambda \subset
\Delta$ and each $\omega \in \Omega$. Let us denote, for each $\Gamma
\subseteq \Delta$, each integer $n \geq \l_\Gamma$ and
$\omega_{-\infty}^{n} \in \Omega_{-\infty}^{n}$,
$$ G_{\Gamma,n}\left(\, \cdot \mid \omega_{\Gamma^{c}} \right)
\;\triangleq\; f_{[l_\Gamma, n]} \Bigl(\, \cdot \,
\ind{\omega_{\Gamma^{c} \cap [l_{\Gamma},n]}} \Bigm|
\omega_{\Gamma_{-}} \Bigr).$$
Definition \reff{comp2}--\reff{comp3} becomes
\equa{\label{lisspe55}
\gamma_{\Delta}\left(h \mid \omega \right) \;=\; \lim_{n \rightarrow
+\infty} \frac{G_{\Delta,n}\left( \ind{\omega_{\Delta}} \mid
\omega_{\Delta^{c}}\right)}{G_{\Delta,n}\left( 1 \mid
\omega_{\Delta^{c}}\right)}.}
Using the reconstruction property
\reff{lis5} of LIS with $l=l_\Delta, \; n=l_{\Lambda}-1$ and $m=n$,
we obtain
$$G_{\Delta,n}( \ind{\omega_{\Delta \setminus \Lambda}} \mid
\omega_{\Delta^{c}}) \;=\; G_{\Delta,l_{\Lambda}-1}(
\ind{\omega_{\Delta \setminus \Lambda}} \mid
\omega_{\Delta^{c}}) \times G_{\Lambda,n}( 1 \mid
\omega_{\Lambda^{c}} )$$
and
$$G_{\Lambda,n}( \ind{\omega_{\Lambda}} \mid \omega_{\Lambda^{c}})
\times G_{\Delta,l_{\Lambda}-1} ( \ind{\omega_{\Delta
\setminus \Lambda }} \mid \omega_{\Delta^{c}} ) \;=\;
G_{\Delta,n}( \ind{\omega_{\Delta}} \mid \omega_{\Delta^{c}})\;.$$
Therefore
\equa{\label{lisspe57} \frac{G_{\Lambda,n} (\ind{\omega_{\Lambda}}
    \mid \omega_{\Lambda^{c}} )}
{G_{\Lambda,n}( 1 \mid \omega_{\Lambda^{c}} )}
\times \frac{G_{\Delta,n}( \ind{\omega_{\Delta \setminus
\Lambda}} \mid \omega_{\Delta^{c}} )}{G_{\Delta,n}( 1 \mid
\omega_{\Delta^{c}} )} \;=\; \frac{G_{\Delta,n} (\ind{\omega_{\Delta}}
    \mid \omega_{\Delta^{c}} )}
{G_{\Delta,n}( 1 \mid \omega_{\Delta^{c}})}\;.}
Identity \reff{lisspe53} follows from \reff{lisspe55} and
\reff{lisspe57}.
\smallskip

We proceed with item {\bf 1) (c)}.  By \reff{lisspe51} and the
summability of the bound $\epsilon_k$ [defined in \reff{lisspe1}],
$F_{\Lambda,n}(\omega_{\Lambda} \mid \cdot \, )$ converges uniformly
to $\gamma_{\Lambda}(\omega_{\Lambda} \mid \cdot \, )$. As each
$F_{\Lambda,n}$ is continuous on $\Omega$, so is
$\gamma_{\Lambda}^{f}$. Let us fix $k_0$ such that $\epsilon_k<1$ for
$k\ge k_0$.  By \reff{lisspe1} and the lower bound in \reff{eq:3}
\[
\gamma_{\Lambda}^{f} \;\ge\; F_{\Lambda,k_0}\, \prod_{k=k_0}^\infty
(1-\epsilon_k)\;.
\]
The right-hand side is strictly positive on $\Omega$ due to the
non-nullness of $F_{\Lambda,k_0}$ and the summability of the
$\epsilon_k$.  Hence $\gamma_{\Lambda}^{f}$ is non-null on $\Omega$.
\smallskip

To prove assertion {\bf 2)(a)} we consider $\mu \in
\mathcal{G}\left(f^{(V,\omega)}\right)$ and denote
$$G_{\Lambda,n}^{(V, \omega)}\left( \, \cdot \mid \sigma_{\Lambda^{c}}
\right) \;\triangleq\; f_{[l_{\Lambda}, n]} \Bigl( \, \cdot \,
\ind{\sigma_{\Lambda^{c} \cap [l_{\Lambda},n]}} \Bigm| \omega_{V_{-}}
\, \sigma_{\Lambda_{-} \setminus V_{-}}\Bigr)$$
for all $\Lambda \in \mathcal{S}_{b} : \Lambda \subset V$ and $\omega,
\sigma : \omega_{V_{-}} \, \sigma_{\Lambda_{-} \setminus V_{-}} \in
\Omega_{\Lambda_{-}}$. By a straightforward extension of
\reff{lisspe55}, the dominated convergence theorem and the consistency
of $\mu$ with respect to $G_{\Lambda,n}^{(V,\omega)}$
\begin{align*}
  \mu \, \gamma_{\Lambda}^{f}\left( \ind{\omega_{\Lambda}}\right) & =
  \lim_{n \rightarrow +\infty} \mu \, G_{\Lambda,n}^{(V,\omega)}
  \left( \frac{G_{\Lambda,n}^{(V,\omega)}\left( \ind{\omega_{\Lambda}}
        \mid \cdot
        \,\right)}{G_{\Lambda,n}^{(V,\omega)}\left( 1 \mid \cdot \, \right)}\right)\\
  & = \lim_{n \rightarrow +\infty} \mu \, G_{\Lambda,n}^{(V,\omega)}
  \left( \ind{\omega_{\Lambda}}\right).
\end{align*}
Applying the consistency hypothesis a second time we obtain
$\mu\,\gamma_{\Lambda}^{f}=\mu$.
\smallskip

Assertion {\bf 2)(b)} is an immediate consequence of 2) (a) and of the
fact that $\left|\mathcal{G}(\gamma)\right|=1$ for all
$\gamma\in\Pi_1$.
\smallskip

Finally we prove {\bf 2) (c)}. Let $f^{1}$ and $f^{2}$ be two LIS on
$(\Omega, \mathcal{F})$, both in $b^{-1}(\Pi_{1})$, and such that
$\gamma^{f^{1}} = \gamma^{f^{2}}.$ By 2) (c), $\mu^{f^{1}} =
\mu^{f^{2}} \triangleq \mu$. The non-nullness of $f^{1}$ and $f^{2}$
on $\Omega$ implies that $\mu$ charges all open sets in $\Omega$.
Therefore, $f_{\Lambda}^{1}$ and $f_{\Lambda}^{2}$ coincide, on
$\Omega$, with the unique continuous realization of $E_{\mu}\left( \,
  \cdot \mid \mathcal{F}_{\Lambda_{-}} \right)$. $\quad \Box$
\subsection{Specification $\rightsquigarrow$ LIS}

Let us introduce the \emph{spread} of a (bounded) function $h$ on $\Omega$:
$$ \Spr{(h)} \;=\; \sup (h) - \inf (h)\;.$$

\begin{lemm}\label{spelislem1}
\begin{fleuve1}
\item Let $\gamma$ be a specification on $\Omega$.
\begin{fleuvea}
\item
If there exists an exhausting sequence of regions
$\Lambda_n\subset\mathbb{Z}$ such that
\equa{\label{spelis3} \lim_{n \rightarrow +\infty} \Spr
\left( \gamma_{\Lambda_{n}} h \right)=0}
for each continuous $\mathcal{F}$-measurable function $h$, then
$\left|\mathcal{G}(\gamma)\right|\le 1$.
\item If $\gamma$ is continuous and
  $\left|\mathcal{G}(\gamma)\right|\le 1$, then \reff{spelis3} holds
  for all exhausting sequences of regions $\Lambda_n\subset\mathbb{Z}$
  and all continuous $\mathcal{F}$-measurable function $h$.
\end{fleuvea}
\item Let $f$ be a LIS on $\Omega$

\begin{fleuvea}
\item If for each $i\in \mathbb{Z}$ and each continuous
$\mathcal{F}_{\leq i}$-measurable continuous function $h$
\equa{\label{r.fl} \lim_{n \rightarrow +\infty} \Spr \left(
f_{[i-n,i]}h\right)=0\;,}
then $\left|\mathcal{G}(f)\right|\le 1$.
\item If $f$ is continuous and $\left|\mathcal{G}(f)\right|\le 1$,
  then \reff{r.fl} is verified for all $i\in \mathbb{Z}$ and all
  continuous $\mathcal{F}_{\leq i}$-measurable continuous function $h$.
\end{fleuvea}
\end{fleuve1}
\end{lemm}

\proof We proof part 1), the proof of 2) is similar.
The obvious spread-reducing relation
$$
\inf_{\widetilde{\omega} \in \Omega} h(\widetilde{\omega}) \;\leq\; (\gamma_{\Lambda} h)(\omega) \;\leq
\;\sup_{\widetilde{\omega} \in
  \Omega} h(\widetilde{\omega})\;,$$
valid for every bounded measurable function $h$ on $\Omega$ and every
configuration $\omega \in \Omega$, plus the consistency condition
\reff{gibbs5} imply that the sequence
$\{\sup(\gamma_{\Lambda_{n}}h)\}$ is decreasing (and bounded below by
$\inf h$), while the sequence $\{ \inf (\gamma_{\Lambda_{n}} h)\}$
is increasing (and bounded above by $\sup h$). Therefore, if
$\mu,\nu\in\mathcal{G}(\gamma)$,
\equa{\label{spelisn7.0}
\mu(h)-\nu(h) \;\le\; \sup(\gamma_{\Lambda_{n}} h)
- \inf(\gamma_{\Lambda_{n}} h)}
for every $n$, which yields
\equa{\label{spelisn7.1}
\left|\mu(h)-\nu(h)\right| \;\le\; \Spr(\gamma_{\Lambda_{n}} h)
}
for every $n$.  This proves item {\bf 1)(a)}.
\smallskip

Regarding {\bf 1)(b)}, we observe that, as $\Omega$ is compact,
  there exist optimizing boundary conditions $\left\{ {\sigma}^{(n)}
  \right\}$ and $\left\{ {\eta}^{(n)} \right\}$ such that
  $(\gamma_{\Lambda_{n}} h) \left({\sigma}^{(n)}\right) = \sup(
  \gamma_{\Lambda_{n}} h).$ and
  $(\gamma_{\Lambda_{n}} h) \left({\eta}^{(n)}\right) = \inf(
  \gamma_{\Lambda_{n}} h).$
(Of course, both sequences of boundary conditions
depend on $h$).  Let $\overline\rho$ and $\underline\rho$ be
respective accumulation point of the sequences of measures
$\left\{ \gamma_{\Lambda_{n}} ( \, \cdot \mid \sigma ^{(n)}) \right\}$
  and
$\left\{ \gamma_{\eta_{n}} ( \, \cdot \mid \sigma ^{(n)}) \right\}$
(they exist by compactness).  Then,
$\overline\rho, \underline\rho\in \mathcal{G}(\gamma)$ (due to the
  continuity of $\gamma$) and
\equa{\label{spelisn8}
\lim_n\Spr(\gamma_{\Lambda_{n}} h) \;\le\;
\overline\rho(h) - \underline\rho(h)\;.
}
Hence the uniqueness of the consistent measure implies
\reff{spelis3}. We learnt this argument from Michael Aizenman (private
communication). $\square$
\medskip

Our last auxiliary result refers to the following notion.
\begin{defi}
A \embf{global specification} $\gamma$ over $(\Omega, \mathcal{F})$ is a family of probability kernels $\left\{
\gamma_{V}\right\}_{V \subset \mathbb{Z}}$, $\gamma_{V} : \mathcal{F} \times \Omega \rightarrow [0,1]$ such that for
all $V \subset \mathbb{Z}$
\begin{fleuvea}
\item For each $\displaystyle{A \in \mathcal{F}, \; \gamma_{V}(A \mid
\cdot \,) \in \mathcal{F}_{V^{c}}}$.
\item For each $\displaystyle{B \in \mathcal{F}_{V^{c}} \text{ and }
\omega \in \Omega, \; \gamma_{V}(B \mid \omega) = \ind{B}(\omega)}$.
\item For each $W \subset \mathbb{Z} : W \supset V, \; \gamma_{W} \,
\gamma_{V} = \gamma_{W}.$
\end{fleuvea}
\end{defi}

\begin{prop}\label{propglospe1}
Let $\gamma$ be a continuous specification over $(\Omega,\mathcal{F})$
which satisfies a HUC. Then $\gamma$ can be extended into a continuous
global specification such that for every subset $V \subset
\mathbb{Z}$,
\equa{\label{propglospe2} \gamma_{V} \left( h \mid
\omega_{V^{c}}\right) \;\triangleq\; \lim_{\Lambda \uparrow V}
\gamma_{\Lambda} \left( h \mid \omega \right)}
for all continuous functions $h \in \mathcal{F}$ and all $\omega \in
\Omega$. Moreover for all $V \subset \mathbb{Z}$ and all $\omega \in
\Omega$,
\equa{\label{propglospe3} \mathcal{G}\left(
\gamma^{(V,\omega)}\right) \;=\; \left\{\gamma_{V} \left( \, \cdot \mid
\omega\right) \right\}.}
\end{prop}
Georgii (1988)\nocite{geo88} gives a proof of this proposition in the
Dobrushin regime (Theorem 8.23). The same proof extends, with minor
adaptations, under a HUC (see Fern\'andez and Pfister,
1997\nocite{fer97}).

\paragraph{Proof of Theorem \protect \ref{th2}}\mbox{}
\smallskip

Items {\bf1) (a)--(b)} are proven in Proposition \reff{propglospe1}.
\smallskip

There are three things to prove regarding {\bf 1) (c)}:
\smallskip\par\noindent \emph{(i) Continuity of} $f^{\gamma}$.  This
is, in fact, an application of Proposition \reff{propglospe1}.
\smallskip\par\noindent \emph{(ii) Non-nullness of} $f^{\gamma}$.
Consider $\Lambda \in \mathcal{S}$, $\omega \in \Omega$, $n \geq
m_{\Lambda}$ and $k \geq 0$. By the non-nullness and the continuity of
$\gamma$ and the compactness of $\Omega_{\Lambda^{c}}$,
there exists $\widetilde{\omega} \in \Omega_{\Lambda^{c}}$ such that
\[
0\;<\; \gamma_{\Lambda}\left(\omega_{\Lambda} \mid
  \widetilde{\omega}\right) \;=\; \inf_{\omega \in \Omega_{\Lambda}^{c}}
\gamma_{\Lambda}\left(\omega_{\Lambda} \mid \omega\right)\;\triangleq\;
c\left(\Lambda,\omega_{\Lambda} \right).
\]
Therefore by the consistency of $\gamma$
\begin{eqnarray*}
f_{\Lambda}^{\gamma}\left( \omega_\Lambda \mid \omega_{\Lambda_{-}}
\right)& =& \lim_{k\to\infty} \gamma_{\left[l_{\Lambda},n+k\right]}
\left(\omega_\Lambda \mid \omega \right)  \;= \;
\lim_{k\to\infty} \Bigl( \gamma_{\left[l_{\Lambda},n+k\right]}
\gamma_\Lambda\Bigr) \left(\omega_\Lambda \mid \omega \right)  \\[5pt]
& \ge & c(\Lambda,\omega_{\Lambda} ) \;>\; 0\;.
\end{eqnarray*}
\smallskip\par\noindent \emph{(iii) Hereditary uniqueness}.  Let us
fix $\omega\in\Omega$ and $V\in\{[j,+\infty[$, $j\in \mathbb{Z}\}\cup \mathbb{Z}$.  For each $i \in \mathbb{Z}$ and $h
\in \mathcal{F}_{\leq i}$. We have
\[
\Spr \left( f_{[i-k,i]}^{\gamma\,(V,\omega)} h \right)
\;\leq\; \lim_{n \rightarrow +\infty} \Spr \left(
\gamma_{[i-k,i+n]}^{(V,\omega)} h \right)\;.
\]
As, by hypothesis, each specification $\gamma^{(V,\omega)}$ admits an
unique Gibbs measure, it follows from lemma \ref{spelislem1} 1) (b) that
\[
\lim_{k \rightarrow +\infty}
\Spr \left( f_{[i-k,i]}^{\gamma\,(V,\omega)} h \right) \;=\; 0\;.
\]
This proves that
$\left| \mathcal{G} \left( f^{\gamma\,(V,\omega)} \right)\right| = 1$
by lemma \ref{spelislem1} 2) (a).
\smallskip

The uniqueness part of assertion {\bf 2) (a)} is contained in the
just proven hereditary uniqueness.  To show that $\mu^{\gamma} \in
\mathcal{G}(f^{\gamma})$, consider $\Lambda \in \mathcal{S}_{b}$ and
$h$ a continuous $\mathcal{F}_{\leq m_\Lambda}$-measurable
function. By the dominated convergence theorem $$ \mu^{\gamma}
f_{\Lambda}^{\gamma}(h) = \lim_{n \rightarrow +\infty} \int \gamma_{
\Lambda \cup \Lambda_{+}^{(n)}} \left( h \mid \xi \right)
\mu^{\gamma}(d \xi).$$ The consistency of $\mu^{\gamma}$
with respect to $\gamma$ implies, hence, that $\mu^{\gamma} \in
\mathcal{G}(f^{\gamma})$.
\smallskip

To prove assertion {\bf 2) (b)}, let $\gamma^{1}$ and $\gamma^{2}$
such that $f^{\gamma^{1}} = f^{\gamma^{2}}$.  By 2) (a)
$\mu^{\gamma^{1}} = \mu^{\gamma^{2}} \triangleq \mu$. The non-nullness
of $\gamma^{1}$ and $\gamma^{2}$ implies that $\mu$ charges all open
sets on $\Omega$. Therefore for each $\Lambda \in \mathcal{S}$,
$\gamma_{\Lambda}^{1}$ and $\gamma_{\Lambda}^{2}$ coincide with the
unique continuous realization of $E_{\mu}\left( \, \cdot \mid
\mathcal{F}_{\Lambda^{c}}\right)$. $\quad \Box$
\medskip

\textbf{Proof of Theorem \protect \ref{th3}}\\ To prove item {\bf
(a)}, let us recall one of the equivalent definitions of the
variational distance between probability measures over
$\left(\Omega_{i},\mathcal{F}_{i} \right)$\\
\[
\left\| \mu - \nu \right\| \;=\;
\sup_{h \in \mathcal{F}_{i}}
\frac{\left| \mu(h) - \nu(h) \right|}{\Spr(h)}\;.
\]
For a proof of this result see for example Georgii (1988) (section
8.1)\nocite{geo88}. By the consistency of
$\stackrel{\circ}{\gamma_{i}}$ with respect to
$\displaystyle{\gamma_{[i,i+k]}, \; k \geq 0}$
\[
\stackrel{\circ}{f_{i}^{\gamma}}\left( \, \cdot \mid
  \omega_{\infty}^{i-1}\right) \;\triangleq\; \lim_{k \rightarrow +\infty}
\stackrel{\circ}{\gamma}_{[i,i+k]}(\, \cdot \mid \omega)\; = \;\lim_{k
  \rightarrow +\infty} \gamma_{[i,i+k]} \,
\stackrel{\circ}{\gamma_{i}} (\, \cdot \mid \omega)\;.
\]
Therefore, by dominated convergence,
 \equa{\label{spelis53} C_{ij}(f^{\gamma}) \;\leq\;
  \sup_{\substack{ \xi, \eta \in \Omega\\[1pt]
\xi_{-\infty}^{i-1}
      \stackrel{\neq j}{=} \eta_{-\infty}^{i-1}}}\left\|
    \stackrel{\circ}{\gamma_{i}}(\, \cdot \mid \xi) -
    \stackrel{\circ}{ \gamma_{i}}(\, \cdot \mid \eta) \right\|.} Since
$\gamma$ is continuous, we can do an infinite telescoping of
\reff{spelis53} to obtain
\[
C_{ij}(f^{\gamma}) \;\leq\; \sum_{\substack{k=j\\ \text{or } k>i}}
C_{ik}(\gamma)\;.
\]
Thus
\[
\sum_{j:j<i} C_{ij}(f^{\gamma}) \;\leq\; \sum_{j:j \neq i}
C_{ij}(\gamma)<1\;.
\]
\smallskip

To show assertion {\bf (b)}, consider $\gamma \in \Pi_{2}$ for which
there exists a constant $K>0$ such that for every cylinder set $A =
\{x_{l}^{m}\} \in \Omega_{l}^{m}$ there exist integers $n,p$
satisfying
\[
 \gamma_{[n,p]}(A \mid \xi) \;\geq\; K \; \gamma_{[n,p]}(A
\mid \eta) \quad \text{for all } \xi, \eta \in \Omega\;.
\]
Hence, by consistency of $\gamma$, we have that for some fixed $\sigma
\in \Omega$ and for each $k \geq 0$
\[
\gamma_{[n,p+k]}(A \mid \xi)  \; =\;
 \int \gamma_{[n,p]}(A \mid \omega)\, \gamma_{[n,p+k]}(d \omega \mid \xi)
\; \geq\; K \, \gamma_{[n,p]}(A \mid \sigma)\;.
\]
In a similar way we obtain
\[
\gamma_{[n,p+k]}(A \mid \eta) \;\leq\; \frac{1}{K} \,
\gamma_{[n,p]}(A \mid \sigma)\;.
\]
We conclude that for each $k\ge 0$
\[
\gamma_{[n,p+k]}(A \mid \xi) \;\geq\; K^{2} \,
\gamma_{[n,p+k]}(A \mid \eta)\;.
\]
Letting $k\to\infty$ we obtain, due to definition \reff{comp5}, that
$f_{[n,m]}^{\gamma}(A \mid \xi) \geq K^{2} f_{[n,m]}^{\gamma}(A \mid
\eta)$. $\quad \Box$
\subsection{LIS $\leftrightsquigarrow$ specification}
\textbf{Proof of Theorem \protect \ref{th5}}\\
Assertion {\bf 1) (a)} is a direct consequence of inequality
\reff{lisspe51} of Lemma \ref{lemproof} with $\Lambda = \{k\}$
and $n=j-1$.
\smallskip

Assertion {\bf 1) (b)} follows from the $n\to\infty$ limit of
inequalities \reff{eq:rr.1} and the fact that $0\le F_{\Lambda,n}\le
1$.
\smallskip

To prove assertion {\bf 2)}, let $k,j \in \mathbb{Z}$ such that $j<k$
and consider $\omega, \sigma \in \Omega$ such that $\omega
\stackrel{\neq j}{=} \sigma$. As a direct consequence of definitions
\ref{def9}--\ref{def13} we have that, for all $i \geq k$,
\equa{\label{exp5}\left(1 -
c_{j}^{\omega}(\gamma_{i})^{-1} \,
\delta_{j}^{\omega}(\gamma_{i})\right) \times \gamma_{i}(\omega_{i}
\mid \omega_{-\infty}^{i-1} \, \omega_{i+1}^{+\infty})
\;\leq\; \gamma_{i}(\sigma_{i} \mid \sigma_{-\infty}^{i-1} \,
\sigma_{i+1}^{+\infty})
}
and
\equa{\label{exp7}
\gamma_{i}(\sigma_{i} \mid \sigma_{-\infty}^{i-1} \,
\sigma_{i+1}^{+\infty})
\;\leq\; \left(1 + c_{j}^{\omega}(\gamma_{i})^{-1} \,
\delta_{j}^{\omega}(\gamma_{i})\right) \times \gamma_{i}(\omega_{i}
\mid \omega_{-\infty}^{i-1} \, \omega_{i+1}^{+\infty})\;.
}
By the
specification reconstruction formula \reff{gibbs37} with $\Lambda =
\{n+1\}$ and $\Gamma=[l_\Lambda,n]$ we have
\[
\gamma_{[l_\Lambda,n+1]}\left( \sigma_\Lambda \mid
\sigma_{\Lambda_-} \, \sigma_{m_\Lambda+1}^{+\infty}\right) \;=\;
\sum_{\sigma_{m_\Lambda+1}^{n}} \frac{\displaystyle{ \gamma_{n+1}\left(
\sigma_{n+1} \mid \sigma_{-\infty}^{n} \,
\sigma_{n+2}^{+\infty}\right)}}{\displaystyle{\sum_{\xi_{n+1}}\frac{\gamma_{n+1}
\left( \xi_{n+1} \mid \sigma_{-\infty}^{n} \, \sigma_{n+2}^{+\infty}
\right)}{\gamma_{[l_\Lambda,n]}\left( \sigma_{l_\Lambda}^{n} \mid
\sigma_{\Lambda_-} \, \xi_{n+1} \,
\sigma_{n+2}^{+\infty}\right)}}}.
\]
Using \reff{exp5} and \reff{exp7}
it is easy to show, by induction over $n \geq m_\Lambda+1$, that
\[
\gamma_{[l_\Lambda,n]}\left( \omega_{\Lambda} \mid \xi_{\Lambda_-} \,
\omega_{n+1}^{+\infty}\right) \;\leq\;
\gamma_{[l_\Lambda,n]}\left( \omega_{\Lambda}
\mid \eta_{\Lambda_-} \, \omega_{n+1}^{+\infty}\right) \times
\prod_{i=k}^{n} \frac{1-c_{j}^{\omega}(\gamma_{k})^{-1} \,
\delta_{j}^{\omega}(\gamma_{i})}{1+c_{j}^{\omega}(\gamma_{k})^{-1} \,
\delta_{j}^{\omega}(\gamma_{i})}
\]
for all $\xi, \eta \in \Omega : \xi \stackrel{\neq j}{=} \eta
\stackrel{\neq j}{=} \omega$. Taking the limit when $n$ tends to
infinity, we obtain 2). $\quad \Box$\\
\newline
\textbf{Proof of Theorem \protect \ref{th4}}\\
For the proof of item {\bf (a)} we consider $\gamma \in \Pi_{2}$ such
that $f^{\gamma} \in \Theta_{1}$ and fix $\Lambda \in \mathcal{S}$ and
$\omega\in\Omega$.  By definition of the maps $b$ and $c$ [see
\reff{comp2}--\reff{comp3} and \reff{comp5}], we have that
\equa{\label{comp16}
\gamma_{\Lambda}^{f^{\gamma}} \left( \omega_{\Lambda} \mid \omega
\right) \;=\;
\lim_{n \rightarrow +\infty} \lim_{k \rightarrow +\infty}
\frac{\displaystyle{ \gamma_{[l_\Lambda,n+k]} \left( \omega_{\Lambda}
\, \omega_{\Lambda^{c} \cap [l_\Lambda,n]} \mid \omega_{\Lambda_{-}} \,
\omega_{n+k+1}^{+\infty} \right)}}
{\displaystyle{ \gamma_{[l_\Lambda,n+k]} \left( \omega_{\Lambda^{c}
\cap [l_\Lambda,n]} \mid \omega_{\Lambda_{-}} \, \omega_{n+k+1}^{+\infty}
\right)}}\;.
}
The consistency of $\gamma_{\Lambda}$ and $\gamma_{[l_\Lambda,n+k]}$
 implies

\alig{\label{comp17} &\gamma_{[l_\Lambda,n+k]} \left( \omega_{\Lambda} \, \omega_{\Lambda^{c} \cap [l_\Lambda,n]} \mid
\omega_{\Lambda_{-}} \, \omega_{n+k+1}^{+\infty}
\right)\nonumber \;=\\[5pt]
& \sum_{\xi_{n+1}^{n+k}} \gamma_{\Lambda} \left( \omega_{\Lambda} \mid \omega_{\Lambda^{c} \cap ]-\infty,n]} \,
\xi_{n+1}^{n+k} \, \omega_{n+k+1}^{+\infty} \right) \, \gamma_{[l_\Lambda,n+k]} \left( \omega_{\Lambda^{c} \cap
[l_\Lambda,n]} \, \xi_{n+1}^{n+k} \mid \omega_{\Lambda_{-}} \, \omega_{n+k+1}^{+\infty}
\right)\;.\nonumber\\
\
}
By continuity of $\gamma_{\Lambda}\left(\omega_{\Lambda} \mid \cdot \,
\right)$ we have that, for each $\varepsilon >0$,
\[
\Bigl| \gamma_{\Lambda} \left( \omega_{\Lambda} \mid
\omega_{\Lambda^{c} \cap ]-\infty,n]} \, \xi_{n+1}^{n+k} \,
\omega_{n+k+1}^{+\infty} \right) - \gamma_{\Lambda}\left(
\omega_{\Lambda} \mid \omega \right) \Bigr|\; <\; \varepsilon
\]
for $n$ large enough uniformly in $k$. Combining this with
\reff{comp16}--\reff{comp17} we conclude that
\[
\left|  \gamma_{\Lambda}^{f^{\gamma}}\left(\omega_{\Lambda} \mid
\omega \right) - \gamma_{\Lambda}\left( \omega_{\Lambda} \mid
\omega \right) \right| \;<\; \varepsilon
\]
for every $\varepsilon >0$.  Therefore $\gamma^{f^{\gamma}} = \gamma$.
\smallskip

To prove item {\bf (b)}, consider $f \in \Theta_{1} $ such that
$\gamma^{f} \in \Pi_{2}$ and fix $\Lambda \in \mathcal{S}_{b}$ and
$\omega \in \Omega$. Let us denote $V = [l_\Lambda,+\infty[$. Since
$\gamma^{f}$ satisfies a HUC equation, \reff{propglospe2} and
definition \reff{comp5} yield
\equa{\label{comp37}
f_{\Lambda}^{\gamma^{f}} \left( \omega_{\Lambda} \mid
\omega_{\Lambda_{-}} \right) \;\triangleq\; \lim_{n \rightarrow +
\infty} \gamma_{[l_\Lambda,m_\Lambda+n]}^{f} \left( \omega_{\Lambda}
\mid \omega \right) \;= \; \gamma_{V}^{f} \left( \omega_{\Lambda} \mid
\omega_{\Lambda_{-}} \right)\;.
}
Combining \reff{propglospe3} with assertion 2) (a) of Theorem
\ref{th1} we obtain that $$\mathcal{G}\left(f^{(V,\omega)} \right) =
\left\{\gamma_{V}^{f}\left(\, \cdot \mid \omega_{\Lambda_{-}}
\right)\right\}.$$ Therefore $$\gamma_{V}^{f}\left(\omega_{\Lambda}
\mid \omega_{\Lambda_{-}}\right) = \gamma_{\Lambda}^{f}\left(
f_{\Lambda}^{(V,\omega)}\left(\omega_{\Lambda} \mid \cdot \, \right)
\mid \omega_{\Lambda_{-}}\right) = f_{\Lambda}\left(\omega_{\Lambda}
\mid \omega_{\Lambda_{-}} \right).$$ The last equality is a
consequence of the definition \reff{lis93}. By \reff{comp37} this
implies that
$$f_{V}^{\gamma^{f}}\left( \omega_{\Lambda} \mid \omega_{\Lambda_{-}}\right) = f_{\Lambda}\left( \omega_{\Lambda} \mid
\omega_{\Lambda_{-}}\right).$$
Item {\bf (c)} is a direct consequence of Theorem \ref{th5} and the
following result. $\Box$

\begin{lemm}\label{explem1}
Let $h : \mathbb{R}^{+} \rightarrow \mathbb{R}^{+}$ be a decreasing
function and $\left(u_{i}\right)_{i \in \mathbb{N}}$ be a sequence
taking values in $]0,1[$ for which there exists $m\geq 0$ such that
$u_{i} \leq m \, h(i)$. Then there exists $M\geq0$ such that
\[
1 - \prod_{i=k}^{+\infty} \frac{1-u_{i}}{1+u_{i}}
\; \leq\; M \, H(k-1)\;,
\]
where $\displaystyle{H(x) = \int_{x}^{+\infty} h(t) \, dt}$.
\end{lemm}
The proof is left to the reader. $\Box$
\appendix
% Appendix : Singleton consistency  for Gibbs measure
\section{Singleton consistency for Gibbs measures}\label{A2}
In this appendix we work in a more general setting than in the paper. We consider a general measurable space
$(E,\mathcal{E})$ (not necessarily finite or even compact) and a subset $\Omega$ of $E^{\mathbb{Z}^{d}}$ for a given
$d\ge 1$. The space $\Omega$ is endowed with the projection $\mathcal{F}$ of the product $\sigma$-algebra associated to
$E^{\mathbb{Z}^{d}}$. We also consider a family of \emph{a priori} measures $\lambda=\left(\lambda^{i}\right)_{i \in
\mathbb{Z}^{d}}$ in $\mathcal{M}\left( E, \mathcal{E}\right)$ and their products $\lambda^{\Lambda} \triangleq
\bigotimes_{i \in \Lambda} \lambda^{i}$ for $\Lambda\subset\mathbb{Z}^d$. We denote by $\left( \lambda_{\Lambda}
\right)_{\Lambda \in \mathcal{S}}$ the family of measure kernels defined over $\left( \Omega, \mathcal{F}\right)$ by
\equa{\label{gibbs24}
\lambda_{\Lambda}(h \mid \omega) =
\left(\lambda^{\Lambda} \otimes \delta_{\omega_{\Lambda^{c}}}
\right)(h)
}
for every measurable function $h$ and configuration
$\omega$. These kernels satisfy the following identities for every
$\Lambda \in \mathcal{S}$:
\equa{\label{gibbs25}
\lambda_{\Lambda}(B \mid \cdot \,) \;=\; \ind{B}(\,
\cdot \,), \; \forall \, B \in \mathcal{F}_{\Lambda^{c}}
}
and
\equa{\label{gibbs26}
\lambda_{\Lambda \cup \Delta} \;=\;
\lambda_{\Lambda} \lambda_{\Delta}\;, \quad \forall \, \Delta \in
\mathcal{S} : \Lambda \cup \Delta = \emptyset.
}
\begin{theo}\label{gibbstheo1}\label{thap}
Let $\lambda$ be as above and $\left( \gamma_{i} \right)_{i \in
\mathbb{Z}^{d}}$ be a family of probability kernels on $\mathcal{F}
\times \Omega_{i}$ such that
\begin{fleuve1}
\item For each $i \in \mathbb{Z}^{d}$ and for some measurable function
$\rho_{i}$,
\equa{\label{gibbs21b}\gamma_{i} \;=\; \rho_{i} \lambda_{i}\;.}
\item The following properties hold:
\begin{fleuvea}
\item Normalization on $\Omega$: for every $i$ in $\mathbb{Z}^{d}$,
\equa{\label{gibbs35}
\left( \lambda_{i} \, (\rho_{i})
\right)(\omega)\;=\;1\;, \quad \forall \, \omega \in \Omega\;.
}
\item Bounded-positivity on $\Omega$: for every $i,j \in
\mathbb{Z}^{d}$,
\equa{\label{gibbs27}
 \inf_{\omega \in \Omega}
\lambda_{j} \left( \rho_{j} \, \rho_{i}^{-1} \right)(\omega)\;>\;0
}
and
\equa{\label{gibbs27.1}
\sup_{\omega \in \Omega} \lambda_{j} \left( \rho_{j} \,
\rho_{i}^{-1} \right)(\omega) \;<\; +\infty\;.
}
\item Order-consistency on $\Omega$: for every $i,j$ in
  $\mathbb{Z}^{d}$ and every $\omega\in\Omega$,
\equa{\label{gibbs31}
\rho_{ij}(\omega) \;=\; \frac{\rho_{i}}{\lambda_{i} \left( \rho_{i} \,
  \rho_{j}^{-1}\right)}(\omega)
\;= \;\frac{\rho_{j}}{ \lambda_{j} \left( \rho_{j} \,
  \rho_{i}^{-1}\right)}(\omega)\;.
}
\end{fleuvea}
\end{fleuve1}
Then there exists a unique family $\rho= \left\{ \rho_{\Lambda}
\right\}_{\Lambda \in \mathcal{S}}$ of positive measurable functions
on $(\Omega,\mathcal{F})$ such that
\begin{fleuvei}
\item $\gamma \triangleq \left\{
\rho_{\Lambda} \lambda_{\Lambda}\right\}_{\Lambda \in \mathcal{S}}$ is
a specification on $(\Omega,\mathcal{F})$ with
$\gamma_{\{i\}}=\gamma_i$ for each $i\in\mathbb{Z}^d$.
\item $\displaystyle{\rho_{\Lambda \cup \Gamma} = \frac{
\rho_{\Lambda}}{\lambda_{\Lambda} \left( \rho_{\Lambda} \,
\rho_{\Gamma}^{-1} \right)}}$, for all $\Lambda, \Gamma \in
\mathcal{S}$ such that $\Gamma \subset \Lambda^{c}$.
\item
%\equa{\label{gibbs39}
$\mathcal{G}(\gamma) =
\left\{ \mu \in \mathcal{P}(\Omega, \mathcal{F}) : \mu \gamma_{i} =
\mu \text{ for all } i \in \mathbb{Z}^{d} \right\}$.
\item For each $\Lambda\in\mathcal{S}$ there exist constants
$C_\Lambda, D_\Lambda >0$ such that
$C_\Lambda\,\rho_k(\omega) \le \rho_\Lambda(\omega) \le
D_\Lambda\,\rho_k(\omega)$ for all $k\in\Lambda$ and all
$\omega\in\Omega)$.
\end{fleuvei}
\end{theo}
\begin{remas1}{3}
\item This theorem is a strengthening of the reconstruction result given by Theorem 1.33 in Georgii
(1988)\nocite{geo88}.  In the latter, the order-consistency condition \reff{gibbs31} is replaced by the requirement
that the singletons come from a pre-existing specification (which the prescription reconstructs). For finite $E$,
Nahapetian and Dachian (2001)\nocite{dacnah01} have presented an alternative approach where \reff{gibbs31} is replaced
by a more detailed pointwise condition. Their non-nullness hypotheses are also different from ours.
\item
Identity (ii) can be used, in fact, to inductively define the family $\rho$ by adding one site at a time.  In fact,
this is what is done in the proof below.  The inequalities (iv) relate the non-nullness properties of $\rho$ to those
of the original family $\{\rho_i\}_{i\in\mathbb{Z}^d}$. \item In the case $E$ countable, $\lambda_i$=counting measure,
the order-consistency requirement \reff{gibbs31} is automatically verified if the singletons are defined through a
measure $\mu$ on $\mathcal{F}$ in the form
\[
\rho_i(\omega) \;=\; \lim_{n\to\infty}
\frac{\mu(\omega_{V_n})}{\mu(\omega_{V_n\setminus{\{i\}}})}
\]
for an exhausting sequence of volumes $\{V_n\}$.  Indeed, a simple
computation shows that the last two terms in \reff{gibbs31} coincide
with
\[
\lim_{n\to\infty}
\frac{\mu(\omega_{V_n})}{\mu(\omega_{V_n\setminus{\{i,j\}}})}\;.
\]
\end{remas1}
\medskip

\proof
In the following all functions are defined on $\Omega$ or on
a projection of $\Omega$ over a subset of $\mathbb{Z}^{d}$.
\smallskip

Initially we define $\rho$ by choosing a total order for
$\mathbb{Z}^{d}$ and prescribing, inductively, that for each $\Lambda
\in \mathcal{S}$ with $|\Lambda| \geq 2$ and each $\omega\in\Omega$
\equa{\label{gibbs45}
\rho_{\Lambda}(\omega) \;=\; \frac{\rho_{k}}{\lambda_{k}\left(
  \rho_{k} \, \rho_{\Lambda_{k}^{*}}^{-1} \right)}(\omega) \;,
}
where $k=\max \Lambda$ and $\Lambda_{k}^{*}=\Lambda \setminus \{k\}$.
For each $\Lambda, \Gamma \in \mathcal{S}$ such that $\Gamma \subset
\Lambda^{c}$, we will prove, by induction over $|\Lambda \cup
\Gamma|$, that the functions so defined satisfy the following
properties:
\begin{fleuveI1}
\item $\displaystyle{ \inf_{\omega \in \Omega} \lambda_{\Lambda}\left(
\rho_{\Lambda} \, \rho_{\Gamma}^{-1} \right)(\omega)>0 \text{ and }
\sup_{\omega \in \Omega} \lambda_{\Lambda}\left( \rho_{\Lambda} \,
\rho_{\Gamma}^{-1} \right)(\omega)<+\infty}.$
\item $\displaystyle{\rho_{\Lambda \cup \Gamma} = \frac{
\rho_{\Lambda}}{\lambda_{\Lambda} \left( \rho_{\Lambda} \,
\rho_{\Gamma}^{-1} \right)}}.$
\item $\displaystyle{
\lambda_{\Lambda}\left(\rho_{\Lambda}\right)=1.}$
\item If $\mu$ is a probability measure on $(\Omega,\mathcal{F})$ such
that $\mu \left(\rho_{i} \, \lambda_{i} \right) = \mu, \; \forall \, i
\in \Lambda$, then $\displaystyle{\mu \left( \rho_{\Lambda} \,
\lambda_{\Lambda} \right) = \mu}.$
\item $\displaystyle{\left( \rho_{\Lambda \cup \Gamma} \,
\lambda_{\Lambda \cup \Gamma}\right) \left( \rho_{i} \,
\lambda_{i}\right) = \rho_{\Lambda \cup \Gamma} \, \lambda_{\Lambda
\cup \Gamma}, \; \forall \, i \in \Lambda \cup \Gamma}.$
\end{fleuveI1}
Let us first comment why these properties imply the theorem.  It is
clear that properties (I3)--(I5), together with the deterministic
character of $\left( \lambda_{\Lambda} \right)$ on
$\mathcal{F}_{\Lambda^c}$ [property \reff{gibbs25}], imply that
$\left( \rho_{\Lambda} \lambda_{\Lambda}\right)_{\Lambda \in
\mathcal{S}}$ verifies assertions (i)--(iv). Furthermore, if
$\widetilde{\gamma}$ is a specification such that
$\widetilde{\gamma}_{\{i\}} = \gamma_{i}$ for all $i \in
\mathbb{Z}^{d}$ then, by consistency, $\widetilde{\gamma}_{\Lambda}(\,
\cdot \mid \omega) \, \gamma_{i} = \widetilde{\gamma}_{\Lambda}(\,
\cdot \mid \omega)$ for every $\Lambda \in \mathcal{S}, \; i \in
\Lambda$ and $\omega \in \Omega$.  Therefore property (I5) implies
that $\widetilde{\gamma}_{\Lambda}(\, \cdot \mid \omega) \, \left(
\rho_{\Lambda} \lambda_{\Lambda}\right) = \widetilde{
\gamma}_{\Lambda}(\, \cdot \mid \omega),$ that is $ \rho_{\Lambda}
\lambda_{\Lambda}(\, \cdot \mid \omega) =
\widetilde{\gamma}_{\Lambda}(\, \cdot \mid \omega)$. So the
construction is unique.
\bigskip

\noindent
{\bf Initial inductive step} The first non-trivial case is when
$|\Lambda \cup \Gamma|=2$. This implies that $|\Lambda|=|\Gamma|=1$
and hence (I1)--(I3) coincide with hypotheses
\reff{gibbs35}--\reff{gibbs31} while (I4) is trivially true. To prove
(I5), assume that $\Lambda=\{i\}$ and $\Gamma=\{j\}$. By
\reff{gibbs26} and \reff{gibbs31}, we have
\[
\left( \rho_{ij} \, \lambda_{ij} \right)\left( \left(\rho_{i} \,
\lambda_{i} \right)(h)\right) \;=\; \lambda_{j} \left[
\left( \frac{\rho_{i} \, \lambda_{i}}{\lambda_{i}\left( \rho_{i} \,
  \rho_{j}^{-1} \right)} \right)\left( \left(\rho_{i}
\, \lambda_{i} \right)(h)\right)\right]\;.
\]
As the factor $\left(\rho_{i} \, \lambda_{i}\right)(h)/ \lambda_{i}\left( \rho_{i} \, \rho_{j}^{-1} \right)$ is
independent of the configuration at $\{i\}$, the remaining integration with respect to the measure $\rho_i\lambda_i$
disappears due to the normalization condition \reff{gibbs35}.  We obtain
\[
\left( \rho_{ij} \, \lambda_{ij} \right)\left( \left(\rho_{i} \,
\lambda_{i} \right)(h)\right) \;=\; \lambda_{j} \left[
\left( \rho_{i} \, \lambda_{i}\right) \left( \frac{h}{ \lambda_{i}
  \left( \rho_{i} \, \rho_{j}^{-1}
\right)}\right)\right]\;=\; \left( \rho_{ij} \, \lambda_{ij}
\right)(h)\;
\]
\medskip

\noindent
{\bf Inductive step}
%Let us denote for every $\Delta$ in $\mathcal{S}$ and $\Lambda$
%subset of $\Delta$, $\Delta_{\Lambda}^{*} \triangleq \Delta \setminus
%\Lambda$.
We suppose the assertions true for $|\Lambda \cup \Gamma|=n, (n\geq 2)$,
and consider $\Lambda, \Gamma$ such that $\Gamma \subset \Lambda^{c}$
and $|\Lambda \cup \Gamma|=n+1$.
\smallskip

{\bf (I1)} Assume first that $|\Gamma|=1$ and let $k = \max
\Lambda$. Combining the definition \reff{gibbs45} and the property
\reff{gibbs26}, we obtain
\equa{\label{eq:rap.1} \lambda_{\Lambda}\left( \rho_{\Lambda} \,
\rho_{\Gamma}^{-1}\right) \;=\; \lambda_{\Lambda_{k}^{*}}\left(
\frac{\lambda_{k} \left(\rho_{k} \, \rho_{\Gamma}^{-1}\right)}{
\lambda_{k} \left( \rho_{k} \,
\rho_{\Lambda_{k}^{*}}^{-1}\right)}\right)\;.
}
If $|\Gamma|\geq 2$ we consider $l \triangleq \max \Gamma$ and apply
the definition \reff{gibbs45} to obtain
\equa{\label{eq:rap.2}
\lambda_{\Lambda}\left( \rho_{\Lambda} \, \rho_{\Gamma}^{-1}\right)
\;=\; \lambda_{\Lambda}\left( \rho_{\Lambda} \, \rho_{l}^{-1} \,
\lambda_{l} \left( \rho_{l} \, \rho_{\Gamma_{l}^{*}}^{-1}\right)
\right)\;.
}
We can now apply the inductive hypothesis (I1) to the right-hand side
of \reff{eq:rap.1} and \reff{eq:rap.2} to prove (I1) at the next
inductive level.
\smallskip

{\bf (I2)} The argument is symmetric in $\Lambda$ and $\Gamma$, so we
can assume without loss that $k=\max(\Lambda\cup\Gamma)$ belongs to
$\Lambda$.  If $\left|\Lambda\right|=1$ (I2) is just the definition
\reff{gibbs45} applied to $\Lambda\cup\Gamma$.  We assume, hence, that
$|\Lambda| \geq 2$ and consider $j \in \Lambda$ such that $j\neq
k$. By the inductive assumption (I2) we have
\equa{\label{gibbs57}
\rho_{\Lambda} \;=\; \frac{\rho_{\Lambda_{
j}^{*}}}{\lambda_{\Lambda_{j}^{*}} \left(
\rho_{\Lambda_{j}^{*}} \, \rho_{j}^{-1}\right)} \;=\;
\frac{\rho_{j}}{\lambda_{j} \left( \rho_{j} \,
\rho_{\Lambda_{j}^{*}}^{-1}\right)}\;.
}
We first combine the rightmost preceding expression with the
factorization property \reff{gibbs26} to write
\equa{\label{gibbs63}
\lambda_{\Lambda} \left(\rho_{\Lambda} \, \rho_{\Gamma}^{-1}\right)
\;=\; \lambda_{\Lambda_{j}^{*}}\left( \frac{\lambda_{j} \left(
\rho_{j} \, \rho_{\Gamma}^{-1} \right)}{\lambda_{j} \left( \rho_{j} \,
\rho_{\Lambda_{j}^{*} }^{-1}\right)}\right)\;.
}
We now apply once more the inductive assumption (I2) in the form
\equa{\label{gibbs65}
\lambda_{j}\left( \rho_{j} \,\rho_{\Gamma}^{-1}\right) \;=\;
 \rho_{\Gamma \cup \{j\}}^{-1} \, \rho_{j}
}
in combination with the rightmost identity in \reff{gibbs57}, to obtain
\equa{\label{gibbs67}
\lambda_{j}\left( \rho_{j} \, \rho_{\Lambda_{j}^{*}}^{-1}\right)
\;=\; \rho_{\Lambda_{j}^{*}}^{-1} \, \rho_{j} \,
\lambda_{\Lambda_{j}^{*}} \left( \rho_{\Lambda_{j}^{*}} \,
\rho_{j}^{-1} \right)\;.
}
From \reff{gibbs63}--\reff{gibbs67} we get
\[
\lambda_{\Lambda} \left(\rho_{\Lambda} \, \rho_{\Gamma}^{-1}\right)
\;=\; \frac{\lambda_{\Lambda_{j}^{*}}\left( \rho_{\Lambda_{j}^{*}}
\, \rho_{\Gamma \cup \{j\}}^{-1} \right)}{
\lambda_{\Lambda_{j}^{*}} \left( \rho_{\Lambda_{j}^{*}} \,
\rho_{j}^{-1} \right)}\;.
\]
We now use this relation together with the first identity in
\reff{gibbs57} to conclude that
\[
\frac{\rho_{\Lambda}}{\lambda_{\Lambda} \left(\rho_{\Lambda} \,
\rho_{\Gamma}^{-1}\right)} \;=\; \frac{\rho_{\Lambda_{j}^{*}}}{
\lambda_{\Lambda_{j}^{*}}\left( \rho_{\Lambda_{j}^{*}}
\, \rho_{\Gamma \cup \{j\}}^{-1} \right)}\;.
\]
We iterate this formula $\left|\Lambda_j^*\right| -1$ times and we
arrive to
\[
\frac{\rho_{\Lambda}}{\lambda_{\Lambda} \left(\rho_{\Lambda} \,
\rho_{\Gamma}^{-1}\right)} \;=\; \frac{\rho_k}{
\lambda_k\left( \rho_k \, \rho_{(\Lambda\cup\Gamma)^*_k}^{-1} \right)}
\]
which is precisely $\rho_{\Lambda\cup\Gamma}$ according to our
definition \reff{gibbs45}.
\smallskip

{\bf (I3)}  We assume that $|\Lambda| \geq 2$, otherwise (I3) is just
the normalization hypothesis \reff{gibbs35}. Let $k = \max \Lambda$.
Definition \ref{gibbs45} and property \ref{gibbs26} yield
\[
\lambda_{\Lambda} \left( \rho_{\Lambda} \right) \;=\;
\lambda_{k} \left( \frac{\lambda_{\Lambda_{k}^{*}}
\left(\rho_{\Lambda_{k}^{*}}\right)}{\lambda_{\Lambda_{k}^{*}} \left(
\rho_{\Lambda_{k}^{*}} \, \rho_{k}^{-1} \right)}\right)
\;=\; \lambda_{k} \left( \frac{1}{\lambda_{\Lambda_{k}^{*}} \left(
\rho_{\Lambda_{k}^{*}} \, \rho_{k}^{-1} \right)}\right)
\]
where the last identity follows from the inductive hypothesis (I3).
But, as in \reff{gibbs65},
\[
\lambda_{\Lambda_{k}^{*}} \left( \rho_{\Lambda_{k}^{*}} \,
\rho_{k}^{-1}\right) \;=\; \lambda_{k} \left( \rho_{k} \,
\rho_{\Lambda_{k}^{*}}^{-1} \right) \, \rho_{\Lambda_{k}^{*}} \,
\rho_{k}^{-1}\;,
\]
therefore
\[
\lambda_{\Lambda} \left( \rho_{\Lambda} \right) \;=\;
\frac{\lambda_{k} \left(\rho_{k} \, \rho_{\Lambda_{k}^{*}}^{-1}
\right)}{\lambda_{k}\left( \rho_{k} \, \rho_{\Lambda_{k}^{*}}^{-1}
\right)} \;=\; 1\;.
\]
\smallskip

{\bf (I4)} To avoid a triviality we assume that $|\Lambda| \geq 2$.
Let $\mu$ be a probability measure on $(\Omega,\mathcal{F})$ such that
$\mu \left( \rho_{i} \, \lambda_{i} \right) = \mu$ for all $i \in
\Lambda$. Consider $k = \max \Lambda$ and a measurable function
$h$. By the factorization property \reff{gibbs26} of $\lambda_\Lambda$
and the definition \reff{gibbs45} of $\rho_{\Lambda}$, we have
\[
\mu \Bigl(\left( \rho_{\Lambda} \lambda_{\Lambda} \right)(h)\Bigr) \;=\;
\mu \left[ \lambda_{\Lambda_{k}^{*}}\left( \biggl(
\frac{\rho_{k}}{\lambda_{k} \left( \rho_{k} \,
\rho_{\Lambda_{k}^{*}}^{-1} \right)} \, \lambda_{k}\biggr) (h) \right)
\right]\;.
\]
By the inductive hypothesis (I4) $\mu$ is consistent with
$\rho_{\Lambda_{k}^{*}} \lambda_{\Lambda_{k}^{*}}$ and with
$\rho_{k} \lambda_{k}$, thus
\[
\mu\Bigr( \left( \rho_{\Lambda} \lambda_{\Lambda} \right)(h)\Bigr)
\;=\; \mu \left[\left( \rho_{k} \, \lambda_{k} \right) \left(
\rho_{\Lambda_{k}^{*}}^{-1}\;  \biggl(
\frac{\rho_{k}}{\lambda_{k} \left( \rho_{k} \,
\rho_{\Lambda_{k}^{*}}^{-1} \right)} \, \lambda_{k}\biggr) (h)
\right) \right]\;.
\]
But, in the right-hand side, the two innermost integrals with respect to
$\lambda_k$ commute with the external one, so we have
\eqna{
\mu\Bigl( \left( \rho_{\Lambda} \lambda_{\Lambda} \right)(h) \Bigr)
&=& \mu \left[\left( \rho_{k} \, \lambda_{k}\right) \left(
\frac{h}{\lambda_{k} \left( \rho_{k} \,
\rho_{\Lambda_{k}^{*}}^{-1} \right)} \quad \lambda_{k} \left( \rho_{k} \,
\rho_{\Lambda_{k}^{*}}^{-1} \right) \right) \right]\nonumber\\
&=& \mu\Bigl(\left(\rho_k\,\lambda_k\right) (h)\Bigr)\;,\nonumber
}
which proves (I4).
\smallskip

{\bf (I5)} Denote $\Delta = \Lambda \cup \Gamma$ and pick $i,j \in
\Delta$, $i \neq j$ and a measurable function $h$. By \reff{gibbs26}
and (I2) we have
\[
\Bigl( \rho_{\Delta} \, \lambda_{\Delta}\Bigr)
\Bigl( \left( \rho_{i} \, \lambda_{i} \right)(h)\Bigr) \;=\;
\lambda_{j} \left[ \frac{\Bigl( \rho_{\Delta_{j}^{*}} \,
\lambda_{\Delta_{j}^{*}}\Bigr)
\Bigl( \left(\rho_{i} \, \lambda_{i} \right)(h)\Bigr)}
{\lambda_{\Delta_{j}^{*}}\left( \rho_{\Delta_{j}^{*}} \,
\rho_{j}^{-1}\right)}\right]\;.
\]
Therefore, applying inductive assumption (I5) we obtain
\eqna{
\Bigr( \rho_{\Delta} \, \lambda_{\Delta}\Bigr)
\Bigl( \left( \rho_{i} \, \lambda_{i} \right)(h)\Bigr) &=&
 \lambda_{j} \left[ \frac{\left(\rho_{\Delta_{j}^{*}} \,
\lambda_{\Delta_{j}^{*}}\right)(h)}{\lambda_{\Delta_{j}^{*}}\left(
\rho_{\Delta_{j}^{*}} \, \rho_{j}^{-1}\right)} \right]
\nonumber\\[5pt]
&=& \left( \rho_{\Delta} \, \lambda_{\Delta} \right)(h)\;.\quad \Box
}

\bibliographystyle{alpha}
%\bibliography{../biblio}
%

%
%
\end{document}